\documentclass[11pt, reqno, psamsfonts]{amsart}
\pdfoutput=1

\usepackage{amssymb}
\usepackage{amsthm}
\usepackage{amsmath}
\usepackage{latexsym}
\usepackage[T1]{fontenc}
\usepackage[utf8]{inputenc}
\usepackage[russian, french, english]{babel}

\usepackage{graphicx}
\usepackage{wrapfig}
\usepackage[justification=centering, labelfont=bf]{caption}
\usepackage{mathtools}
\usepackage{tikz}
\usepackage{amsbsy}
\usepackage[inline]{enumitem}
\usepackage{mathrsfs}
\usetikzlibrary{shapes,snakes}
\usetikzlibrary{arrows.meta}
\usetikzlibrary{decorations.pathmorphing}
\usetikzlibrary{patterns}
\usepackage{float}
\usepackage{array}
\usepackage{multicol}
\usepackage{stmaryrd}
\usepackage{cancel} %
\usepackage{lmodern}
\usepackage{mathabx}
\usepackage{upgreek}
\usepackage{titlesec}
\usepackage[spacing=true,kerning=true,babel=true,tracking=true]{microtype}
\usepackage[shortcuts]{extdash}
\usepackage[foot]{amsaddr}
\usepackage[left=1in,right=1in,top=1in,bottom=1in,bindingoffset=0cm]{geometry}
\usepackage{bm}
\usepackage{centernot}
\usepackage{mdframed}
\usepackage[hidelinks]{hyperref}
\newcommand{\diff}{\,\mathrm{d}}

\usepackage[backend=biber, style=alphabetic, sorting=nyt, maxnames=100,backref=true]{biblatex}
\title{\sffamily Weak degeneracy of graphs}
\date{}

\author{Anton~Bernshteyn}
\address[Anton Bernshteyn]{\normalfont School of Mathematics, Georgia Institute of Technology, Atlanta, GA, USA}
\email{bahtoh@gatech.edu}

\author{Eugene Lee}
\address[Eugene Lee]{\normalfont Department of Mathematical Sciences, Carnegie Mellon University, Pittsburgh, PA, USA}
\email{eleehuaj@andrew.cmu.edu}

\thanks{Research of the first named author is partially supported by the NSF grant DMS-2045412.}

\newtheoremstyle{bfnote}%
{}{}%
{\slshape}{}%
{\bfseries}{\bfseries.}%
{ }%
{\thmname{#1}\thmnumber{ #2}\thmnote{ \ep{\normalfont{}#3}}}

\theoremstyle{bfnote}
\newtheorem{theo}[equation]{Theorem}
\newtheorem*{theo*}{Theorem}
\newtheorem{prop}[equation]{Proposition}
\newtheorem{lemma}[equation]{Lemma}
\newtheorem{claim}[equation]{Claim}

\newtheorem{conj}[equation]{Conjecture}

\newtheorem*{corl*}{Corollary}

\theoremstyle{definition}
\newtheorem{defn}[equation]{Definition}
\newtheorem*{defn*}{Definition}

\newtheorem*{exmp*}{Example}

\theoremstyle{remark}
\newtheorem*{ques*}{Question}
\newtheorem*{remk*}{Remark}

\makeatletter
\newcommand{\neutralize}[1]{\expandafter\let\csname c@#1\endcsname\count@}
\makeatother

\newenvironment{propcopy}[1]
{\renewcommand{\theequation}{\ref{#1}}%
	\neutralize{equation}\phantomsection
	\begin{prop}}
	{\end{prop}}

\newenvironment{theocopy}[1]
{\renewcommand{\theequation}{\ref{#1}}%
	\neutralize{equation}\phantomsection
	\begin{theo}}
	{\end{theo}}

\newenvironment{scproof}[1][]{\begin{proof}[\textsc{\upshape{Proof}}#1]}{\end{proof}}

\newcommand{\0}{\varnothing}
\newcommand{\set}[1]{\{#1\}}
\newcommand{\N}{{\mathbb{N}}}
\newcommand{\Z}{\mathbb{Z}}

\newcommand{\R}{\mathbb{R}}

\renewcommand{\P}{\mathbb{P}}
\newcommand{\E}{\mathbb{E}}
\renewcommand{\epsilon}{\varepsilon}

\renewcommand{\phi}{\varphi}
\renewcommand{\theta}{\vartheta}
\renewcommand{\leq}{\leqslant}
\renewcommand{\geq}{\geqslant}
\newcommand{\defeq}{\coloneqq}

\newcommand{\bemph}[1]{{\normalfont#1}} %
\newcommand{\ep}[1]{\bemph{(}#1\bemph{)}} %

\newcommand{\emphd}[1]{{\fontseries{b}\selectfont\textsf{#1}}}

\newcommand{\del}{\textsc{Delete}\xspace}
\newcommand{\delsave}{\textsc{DelSave}\xspace}
\newcommand{\dg}{\mathsf{d}}
\newcommand{\wdg}{\mathsf{wd}}
\newcommand{\mad}{\mathrm{mad}}
\newcommand{\ad}{\mathrm{ad}}
\newcommand{\abs}[1]{\left|#1\right|}
\newcommand{\pto}{\dashrightarrow}
\newcommand{\rest}[2]{{{#1}\vert_{#2}}}
\numberwithin{equation}{section}

\usepackage{xspace}
\newcommand{\Lister}{\texttt{Lister}\xspace}
\newcommand{\Painter}{\texttt{Painter}\xspace}

\titleformat{\section}[block]{\large\bfseries\sffamily}{\thesection.}{1ex}{}
\titleformat{\subsection}[block]{\bfseries\sffamily}{\thesubsection.}{1ex}{}
\titleformat{\subsubsection}[runin]{\itshape}{\bfseries\upshape\thesubsubsection.}{1ex}{}[.---]

\titlespacing*{\section}{0pt}{*3}{*1}
\titlespacing*{\subsection}{0pt}{*3}{*1}
\titlespacing*{\subsubsection}{0pt}{*1.5}{*0}

\renewbibmacro{in:}{}

\renewbibmacro*{volume+number+eid}{%
	\printfield{volume}%
	\setunit*{\addnbspace}%
	\printfield{number}%
	\setunit{\addcomma\space}%
	\printfield{eid}}

\DeclareFieldFormat[article]{volume}{\textbf{#1}\space}
\DeclareFieldFormat[article]{number}{\mkbibparens{#1}}

\DeclareFieldFormat{journaltitle}{#1,}
\DeclareFieldFormat[thesis]{title}{\mkbibemph{#1}\addperiod}
\DeclareFieldFormat[article, unpublished, thesis]{title}{\mkbibemph{#1},}
\DeclareFieldFormat[book]{title}{\mkbibemph{#1}\addperiod}
\DeclareFieldFormat[unpublished]{howpublished}{#1, }

\DeclareFieldFormat{pages}{#1}

\DeclareFieldFormat[article]{series}{Ser.~#1\addcomma}

\setlist{topsep=3pt,itemsep=3pt}

\begin{document}
	
	\maketitle
	
	\begin{abstract}
		Motivated by the study of greedy algorithms for graph coloring, we introduce a new graph parameter, which we call \emph{weak degeneracy}. By definition, every $d$-degenerate graph is also weakly $d$-degenerate. On the other hand, if $G$ is weakly $d$-degenerate, then $\chi(G) \leq d + 1$ \ep{and, moreover, the same bound holds for the list-chromatic and even the DP-chromatic number of $G$}. It turns out that several upper bounds in graph coloring theory can be phrased in terms of weak degeneracy. For example, we show that planar graphs are weakly $4$-degenerate, which implies Thomassen's famous theorem that planar graphs are $5$-list-colorable. We also prove a version of Brooks's theorem for weak degeneracy: a connected graph $G$ of maximum degree $d \geq 3$ is weakly $(d-1)$-degenerate unless $G \cong K_{d + 1}$. \ep{By contrast, all $d$-regular graphs have degeneracy $d$.} We actually prove an even stronger result, namely that for every $d \geq 3$, there is $\epsilon > 0$ such that if $G$ is a graph of weak degeneracy at least $d$, then either $G$ contains a $(d+1)$-clique or the maximum average degree of $G$ is at least $d + \epsilon$. Finally, we show that graphs of maximum degree $d$ and either of girth at least $5$ or of bounded chromatic number are weakly $(d - \Omega(\sqrt{d}))$-degenerate, which is best possible up to the value of the implied constant. %
	\end{abstract}
	
	\section{Introduction}

	All graphs in this paper are finite and simple. Recall that for a graph $G$, $\chi(G)$ denotes its \emphd{chromatic number}, i.e., the minimum number of colors necessary to color the vertices of $G$ so that adjacent vertices are colored differently. A well-studied generalization of graph coloring is \emphd{list coloring}, which was introduced independently by Vizing \cite{Viz} and Erd\H{o}s, Rubin, and Taylor \cite{ERT}. In the setting of list coloring, each vertex $u \in V(G)$ is given a set $L(u)$, called its \emphd{list of available colors}. A \emphd{proper $L$-coloring} is then a function $\phi$ defined on $V(G)$ such that:
	\begin{itemize}
		\item $\phi(u) \in L(u)$ for all $u \in V(G)$; and
		\item $\phi(u) \neq \phi(v)$ for all $uv \in E(G)$.
	\end{itemize}
	The \emphd{list-chromatic number} of $G$, denoted by $\chi_\ell(G)$, is the minimum $k$ such that $G$ admits a proper $L$-coloring whenever $|L(u)| \geq k$ for all $u \in V(G)$. Clearly, $\chi_\ell(G) \geq \chi(G)$ for all graphs $G$.
	
	A further generalization of list coloring is \emphd{DP-coloring} \ep{also known as \emphd{correspondence coloring}}, which was recently introduced by Dvo\v{r}\'ak and Postle \cite{DP}. A related notion of \emph{local conflict coloring} was studied independently from the algorithmic standpoint by Fraigniaud, Heinrich, and Kosowski \cite{FHK}. Just as in list coloring, we assume that every vertex $u \in V(G)$ of a graph $G$ is given a list $L(u)$ of colors to choose from. In contrast to list coloring though, the identifications between the colors in the lists are allowed to vary from edge to edge. That is, each edge $uv \in E(G)$ is assigned a matching $C_{uv}$ \ep{not necessarily perfect and possibly empty} from $L(u)$ to $L(v)$. If $\alpha \beta \in C_{uv}$, we say that $\alpha$ \emphd{corresponds} to $\beta$ \ep{under the correspondence $C$}. A \emphd{proper $(L, C)$-coloring} of $G$ is a function $\phi$ defined on $V(G)$ such that:
	\begin{itemize}
		\item $\phi(u) \in L(u)$ for all $u \in V(G)$; and
		\item $\phi(u)\phi(v) \not \in C_{uv}$ for all $uv \in E(G)$.
	\end{itemize}
	List coloring is a special case of this framework where $\alpha \in L(u)$ corresponds to $\beta \in L(v)$ if and only if $\alpha = \beta$. The \emphd{DP-chromatic number} of $G$, denoted by $\chi_{DP}(G)$, is the minimum $k$ such that $G$ admits a proper $(L,C)$-coloring whenever $|L(u)| \geq k$ for all $u \in V(G)$. Again, it is clear from the definition that $\chi_{DP}(G) \geq \chi_\ell(G)$.
	
	In this paper we are interested in greedy algorithms for graph coloring. The basic greedy algorithm considers the vertices of $G$ one at a time. When we get to consider a vertex $u$, we assign to it an arbitrary color, say $\alpha$, from $L(u)$. At this point, to ensure that the coloring is proper, we have to remove the colors corresponding to $\alpha$ from the lists of colors available to the neighbors of $u$. Thus, the list size for every neighbor of $u$ may decrease by $1$, while all the other lists remain unchanged. If throughout this process no list size reduces to $0$ \ep{i.e., if every uncolored vertex always has at least one available color}, then we successfully obtain a proper \ep{DP-}coloring of $G$. This idea is formally captured in the notion of \emph{graph degeneracy}:
	
	\begin{defn}[\textls{Degeneracy}]
		Let $G$ be a graph and let $f \colon V(G) \to \N$ be a function.\footnote{In this paper $\N = \set{0,1,2,\ldots}$ denotes the set of all non-negative integers.} For a vertex $u \in V(G)$, the operation $\del(G, f, u)$ outputs the graph $G' \defeq G - u$ and the function $f' \colon V(G') \to \Z$ given by the formula
		\[
		f'(v) \,\defeq\, \begin{cases}
			f(v) - 1 &\text{if } uv \in E(G);\\
			f(v) &\text{otherwise}.
		\end{cases}
		\]
		An application of the operation $\del$ is \emphd{legal} if the resulting function $f'$ is non-negative, i.e., if $f'(v) \geq 0$ for all $v \in V(G')$. A graph $G$ is \emphd{$f$-degenerate} if it is possible to remove all vertices from $G$ by a sequence of legal applications of the operation $\del$. Given $d \in \N$, we say that $G$ is \emphd{$d$-degenerate} if it is degenerate with respect to the constant $d$ function. The \emphd{degeneracy} of $G$, denoted by $\dg(G)$, is the minimum $d$ such that $G$ is $d$-degenerate.
	\end{defn}
	
	It follows from the above discussion that $\chi_{DP}(G) \leq \dg(G) + 1$ for every graph $G$; because of this, the quantity $\dg(G) + 1$ is sometimes referred to as the \emph{coloring number} of $G$ \cite{EH}. It is not hard to see that a graph $G$ is $d$-degenerate if and only if every nonempty subgraph of $G$ has a vertex of degree at most $d$ \cite[Proposition 5.2.2]{Die}.
	
	The upper bound $\chi_{DP}(G) \leq \dg(G) + 1$ is usually not sharp. For instance, if $G$ is a $d$-regular graph, then $\dg(G) = d$, which implies that $\chi_{DP}(G) \leq d + 1$. However, the only connected $d$-regular graphs $G$ with $\chi_{DP}(G) = d + 1$ are the complete graph $K_{d+1}$ and---if $d=2$---cycles \cite{BKP}. \ep{A curious distinction between DP-coloring and list coloring is that $\chi_\ell(C_n)$ is $2$ if $n$ is even and $3$ if $n$ is odd, while $\chi_{DP}(C_n) = 3$ for all $n \geq 3$ \cite[\S1.1]{DP}.} It is therefore interesting to see if we can modify the greedy coloring procedure to ``save'' some of the colors and get a better bound on $\chi_{DP}(G)$. Here we investigate a particularly simple \ep{but, as we shall see, already quite powerful} way of doing so. %
	
	To motivate our main definition, consider a vertex $u \in V(G)$ and let $w$ be its neighbor. In general, if we assign a color to $u$, then $w$ may lose one of its colors. However, suppose that $|L(u)| > |L(w)|$, i.e., that $u$ has \emph{strictly more} available colors than $w$. In this case, there must be a color in $L(u)$ that does not correspond to any color in $L(w)$, and assigning such a color to $u$ does not affect $L(w)$ \ep{of course, the other neighbors of $u$ may still lose a color}. In this way, we ``save'' an extra color for $w$. This idea naturally leads to the notion that we call \emph{weak degeneracy}:
	
	\begin{defn}[\textls{Weak degeneracy}]\label{defn:weak_deg}
		Let $G$ be a graph and let $f \colon V(G) \to \N$ be a function. For a pair of adjacent vertices $u$, $w \in V(G)$, the operation $\delsave(G, f, u, w)$ outputs the graph $G' \defeq G - u$ and the function $f' \colon V(G') \to \Z$ given by the formula
		\[
		f'(v) \,\defeq\, \begin{cases}
			f(v) - 1 &\text{if } uv \in E(G) \text{ and } v \neq w;\\
			f(v) &\text{otherwise}.
		\end{cases}
		\]
		An application of the operation $\delsave$ is \emphd{legal} if $f(u) > f(w)$ and the resulting function $f'$ is non-negative. A graph $G$ is \emphd{weakly $f$-degenerate} if it is possible to remove all vertices from $G$ by a sequence of legal applications of the operations $\del$ and $\delsave$. Given $d \in \N$, we say that $G$ is \emphd{weakly $d$-degenerate} if it is weakly degenerate with respect to the constant $d$ function. The \emphd{weak degeneracy} of $G$, denoted by $\wdg(G)$, is the minimum $d$ such that $G$ is weakly $d$-degenerate.
	\end{defn}
	
	Again, the above discussion shows that $\chi_{DP}(G) \leq \wdg(G) + 1$ for every graph $G$. Actually, the same bound holds even for the on-line version of DP-chromatic number called \emph{DP-paint number}, which was introduced by Kim, Kostochka, Li, and Zhu \cite{Paint} \ep{see \S\ref{sec:paint} for the definition}:
	
	\begin{prop}\label{prop:paint}
		For every graph $G$,
		\[
		\chi(G) \,\leq\, \chi_\ell(G) \,\leq\, \chi_{DP}(G) \,\leq\, \chi_{DPP}(G) \,\leq\, \wdg(G) + 1,
		\]
		where $\chi_{DPP}(G)$ is the DP-paint number of $G$.
	\end{prop}
	
	It turns out that the simple way of ``saving'' colors using the $\delsave$ operation is sufficient for several non-trivial upper bounds. For example, consider the case of planar graphs. It follows from Euler's formula that planar graphs are $5$-degenerate, which gives a simple proof of their $6$-colorability \ep{and even $6$-DP-colorability}. On the other hand, Thomassen \cite{Thom} proved that every planar graph is $5$-list-colorable, and this result was extended to DP-coloring by Dvo\v{r}\'ak and Postle \cite{DP}. The value $5$ here is optimal as Voigt \cite{Voigt} constructed planar graphs of list-chromatic number exactly $5$. While degeneracy is not sufficient to establish Thomassen's theorem, we show in \S\ref{sec:planar} that weak degeneracy is:
	
	\begin{theo}\label{theo:planar}
		Every planar graph is weakly $4$-degenerate.
	\end{theo}
	
	Next we consider Brooks-type theorems for weak degeneracy. As mentioned earlier, if $d \geq 3$, then the only connected graph $G$ of maximum degree $d \geq 3$ with $\chi_{DP}(G) = d + 1$ is the complete graph $K_{d+1}$. We show that there is a corresponding bound on weak degeneracy:
	
	\begin{theo}\label{theo:Brooks}
		If $G$ is a connected graph of maximum degree $d \geq 3$, then either $G \cong K_{d+1}$ or $G$ is weakly $(d-1)$-degenerate. 
	\end{theo}
	
	More generally, suppose that $G$ is a connected graph and $|L(u)| \geq \deg_G(u)$ for every vertex $u \in V(G)$ \ep{that is, the lower bound on the list size varies depending on the degree of the vertex}. In the list-coloring framework, Borodin \cite{Bor} and, independently, Erd\H{o}s, Rubin, and Taylor \cite{ERT} showed that $G$ is $L$-colorable unless it is a \emphd{Gallai tree}, i.e., a connected graph in which every block is either a clique or an odd cycle. In the DP-coloring setting the same result holds, except that the graphs that need to be excluded are the \emphd{GDP trees}, i.e., connected graphs in which every block is either a clique or a cycle \ep{not necessarily odd} \cite{BKP}. We again establish the corresponding result for weak degeneracy:
	
	\begin{theo}\label{theo:GDP}
		Let $G$ be a connected graph. The following statements are equivalent:
		\begin{enumerate}[label=\ep{\normalfont\arabic*}]
			\item $G$ is weakly $f$-degenerate, where $f(u) = \deg_G(u) - 1$ for all $u \in V(G)$;
			\item $G$ is not a GDP-tree.
		\end{enumerate}
	\end{theo}
	
	Theorem~\ref{theo:Brooks} is an immediate corollary of Theorem~\ref{theo:GDP}. We prove Theorem~\ref{theo:GDP} in \S\ref{sec:GDP}.
	
	Recall that the \emphd{average degree} of a nonempty graph $G$, denoted by $\ad(G)$, is the average of the degrees of the vertices of $G$. Equivalently, we have $\ad(G) = 2|E(G)|/|V(G)|$. The \emphd{maximum average degree} of $G$, denoted by $\mad(G)$, is defined by $\mad(G) \defeq \max_H \ad(H)$, where the maximum is taken over all nonempty subgraphs $H$ of $G$. The maximum average degree of a graph is a natural measure of its local density. There is a close relationship between a graph's maximum average degree and its degeneracy; namely, we have
	\[
	2\dg(G) \,\geq\, \mad(G) \,\geq\, \dg(G).
	\]
	For $d$-regular graphs $G$, $\mad(G) = \dg(G) = d$. By contrast, we show that if $\wdg(G) \geq 3$ and $G$ contains no $(\wdg(G) + 1)$-clique, then $\mad(G) \geq \wdg(G) + \epsilon$, where $\epsilon > 0$ only depends on $\wdg(G)$:
	
	\begin{theo}\label{theo:mad}
		Let $G$ be a nonempty graph. If the weak degeneracy of $G$ is at least $d \geq 3$, then either $G$ contains a $(d+1)$-clique or
		\[
		\mathrm{mad}(G) \,\geq\, d \,+\, \frac{d-2}{d^2 + 2d - 2}. %
		\]
	\end{theo}
	
	Note that Theorem~\ref{theo:mad} is a strengthening of Theorem~\ref{theo:Brooks}, since $\mad(G)$ is at most the maximum degree of $G$. Our proof of Theorem~\ref{theo:mad}, which we present in \S\ref{sec:mad}, relies on Theorem~\ref{theo:GDP} and follows an approach similar to the one used by Gallai \cite{Gallai} to establish a lower bound on the average degree of critical graphs.
	
	As far as lower bounds on weak degeneracy are concerned, a fairly straightforward double counting argument gives the following:
	
	\begin{prop}\label{prop:lower_bound}
		Let $G$ be a $d$-regular graph with $n \geq 2$ vertices. Then $\wdg(G) \geq d - \sqrt{2n}$. %
	\end{prop}
	
	In particular, if $n= O(d)$, then $\wdg(G) \geq d - O(\sqrt{d})$. For example, Proposition~\ref{prop:lower_bound} yields the bound $\wdg(K_{d,d}) \geq d - 2\sqrt{d}$ for $d \geq 2$. Actually, this can be improved to $d - \sqrt{2d} - 1$:%
	
	\begin{prop}\label{prop:Kdd}
		If $G$ is a triangle-free $d$-regular graph with $n \geq 4$ vertices, then $\wdg(G) > d - \sqrt{n} - 1$. In particular, %
		the complete bipartite graph $K_{d,d}$ with $d \geq 2$ satisfies $\wdg(K_{d,d}) > d - \sqrt{2d} - 1$.
	\end{prop}
	
	This should be contrasted with the fact that $\chi(K_{d,d}) = 2$, $\chi_\ell(K_{d,d}) = (1 + o(1))\log_2 d$ \cite{ERT}, and $\chi_{DP}(K_{d,d}) = \Theta(d/\log d)$ \cite{Ber16}. We prove Propositions \ref{prop:lower_bound} and \ref{prop:Kdd} in \S\ref{sec:Kdd}.
	
	It seems plausible that \emph{every} $d$-regular graph has weak degeneracy at least $d - O(\sqrt{d})$; we leave verifying or refuting this supposition as an open problem:
	
	\begin{conj}\label{conj:reg}
		Every $d$-regular graph $G$ satisfies $\wdg(G) \geq d - O(\sqrt{d})$.
	\end{conj}
	
	In view of the above lower bounds, it makes sense to ask, for what classes of graphs $G$ does the upper bound $\wdg(G) \leq d - \Omega(\sqrt{d})$ hold, where $d$ is the maximum degree of $G$? Along these lines, we establish the following results:
	
	\begin{theo}\label{theo:chrom}
		For each integer $k \geq 1$, there exist $c > 0$ and $d_0 \in \N$ such that if $G$ is a graph of maximum degree $d \geq d_0$ with $\chi(G) \leq k$, then $\mathsf{wd}(G) \leq d - c\sqrt{d}$.
	\end{theo}
	
	\begin{theo}\label{theo:girth5}
		There exist $c > 0$ and $d_0 \in \N$ such that if $G$ is a graph of maximum degree $d \geq d_0$ and girth at least $5$, then $\mathsf{wd}(G) \leq d - c \sqrt{d}$.
	\end{theo}
	
	Theorems~\ref{theo:chrom} and \ref{theo:girth5} are proved in \S\ref{sec:chrom} using probabilistic arguments.
	
	We finish the introduction with another conjecture that implies both Theorems~\ref{theo:chrom} and \ref{theo:girth5}:
	
	\begin{conj}\label{conj:clique}
		For each integer $k \geq 1$, there exist $c > 0$ and $d_0 \in \N$ such that if $G$ is a graph of maximum degree $d \geq d_0$ and without a $k$-clique, then $\mathsf{wd}(G) \leq d - c \sqrt{d}$.
	\end{conj}
	
	At present, we do not even know if Conjecture~\ref{conj:clique} holds for $k = 3$, i.e., whether $\wdg(G) \leq d - \Omega(\sqrt{d})$ for triangle-free graphs of maximum degree $d$.
	
	\subsubsection*{Acknowledgments} We are very grateful to the anonymous referees for their comments.
	
	\section{Preliminary results}
	
	In this section we establish several basic results about weak degeneracy that will be used throughout the rest of this article.
	
	\begin{lemma}[\textls{Weak degeneracy is monotone}]\label{lemma:monotone}
		Let $G$ be a weakly $f$-degenerate graph. If $g \colon V(G) \to \N$ is a function such that $g(u) \geq f(u)$ for all $u \in V(G)$, then $G$ is weakly $g$-degenerate.
	\end{lemma}
	\begin{scproof}
		We wish to show that $G$ is weakly $g$-degenerate by removing its vertices via the same sequence of operations that witnesses that $G$ is weakly $f$-degenerate. The only possible issue is that an application of \delsave may become illegal when $f$ is replaced by $g$. Namely, it can happen that $\delsave(G,f,u,w)$ is legal, while $\delsave(G,g,u,w)$ is not due to the fact that $g(u) \leq g(w)$. However, since $f(u) > f(w)$, this means that $g(w) > f(w)$, so instead of $\delsave(G,g,u,w)$ we can simply use $\del(G,g,u)$: this replaces $g(w)$ by $g(w)-1$, which is still at least $f(w)$.
	\end{scproof}
	
	\begin{lemma}[\textls{Weak degeneracy and DP-coloring}]\label{lemma:DP}
		Let $G$ be a weakly $f$-degenerate graph. Suppose that every vertex $u \in V(G)$ is given a list $L(u)$ of available colors and that for each edge $uv \in E(G)$, there is a matching $C_{uv}$ from $L(u)$ to $L(v)$. If $|L(u)| \geq f(u) + 1$ for all $u \in V(G)$, then $G$ admits a proper $(L,C)$-coloring.
	\end{lemma}
	\begin{scproof}
		This statement was essentially established in the introduction \ep{in the discussion preceding Definition~\ref{defn:weak_deg}}. We give a more detailed proof here for completeness. Since $G$ is weakly $f$-degenerate, it is possible to remove all vertices from $G$ by a sequence of legal applications of the operations \del and \delsave. Fix any such sequence $\mathcal{S} = (\mathcal{O}_0, \ldots, \mathcal{O}_{n-1})$. Set $(G_0, f_0) \defeq (G,f)$ and for each $0 \leq i \leq n-1$, let $(G_{i+1},f_{i+1})$ be the result of applying the operation $\mathcal{O}_i$ to $(G_i, f_i)$. We color the vertices of $G$ one by one, in the order in which they are removed by $\mathcal{S}$. Each time a vertex $u$ is assigned a color $\alpha$, we remove the colors corresponding to $\alpha$ from the lists of colors available to the neighbors of $u$, thus ensuring that the resulting coloring is proper. Let $L_i(u)$ be the list of colors available to a vertex $u \in V(G_i)$ at the start of step $i$ \ep{in particular, $L_0(u) \defeq L(u)$}. Throughout our coloring procedure, we will maintain the following property:
		
		\smallskip
		
		\begin{quote}
			\ep{$\mathrm{P}_i$}\quad\quad $|L_i(u)| \geq f_i(u) + 1$ for all $u \in V(G_i)$. %
		\end{quote}
		
		\smallskip
		
		\noindent If we can achieve this, then we will successfully color the entire graph, since no uncolored vertex will ever run out of available colors. Now, property \ep{$\mathrm{P}_0$} holds by assumption. On step $i$, we assume that \ep{$\mathrm{P}_i$} holds and consider two cases.
		
		\smallskip
		
		\textbf{Case 1:} $\mathcal{O}_i = \del(G_i,f_i,u_i)$.
		
		\smallskip
		
		In this case we assign to $u_i$ an arbitrary available color. It is clear that property \ep{$\mathrm{P}_{i+1}$} holds regardless of what color is assigned to $u_i$. %
		
		\smallskip
		
		\textbf{Case 2:} $\mathcal{O}_i = \delsave(G_i,f_i,u_i,w_i)$.
		
		\smallskip
		
		If $|L(w_i)| > f_i(w_i) + 1$, %
		we can, as in Case 1, assign an arbitrary available color to $u_i$. Now suppose that $|L(w_i)| = f_i(w_i) + 1$. %
		Then, by \ep{$\mathrm{P}_i$} and since this application of $\delsave$ is legal, we have %
		$|L_i(u_i)| \geq f_i(u_i) + 1 > f_i(w_i) + 1 = |L_i(w_i)|$. %
		This means that $u_i$ must have an available color $\alpha_i \in L_i(u_i)$ that does not correspond to any color in $L_i(w_i)$. %
		If we assign $\alpha_i$ to $u_i$, then the list of available colors for $w_i$ will not change, and thus \ep{$\mathrm{P}_{i+1}$} will not be violated, as desired.
	\end{scproof}
	
	\begin{lemma}[\textls{Partitioning lemma}]\label{lemma:part}
		Let $G$ be a weakly $f$-degenerate graph. Suppose that functions $f_1$, $f_2 \colon V(G) \to \Z$ satisfy $f_1(u) + f_2(u) = f(u) - 1$ for all $u \in V(G)$. Then there is a partition $V(G) = V_1 \sqcup V_2$ such that the subgraph $G[V_i]$ is weakly $f_i$-degenerate for each $i \in \set{1,2}$.
	\end{lemma}
	\begin{scproof}
		The proof is by induction on $|V(G)|$. If $V(G) = \0$, the statement holds vacuously. Now suppose that $|V(G)| \geq 1$ and the claim holds for all graphs with $|V(G)| - 1$ vertices. Since $G$ is weakly $f$-degenerate, there is a legal application of an operation $\mathcal{O} \in \set{{\del}, {\delsave}}$ that produces a pair $(G',f')$ in which the graph $G'$ is weakly $f'$-degenerate. We consider the two cases depending on whether $\mathcal{O}$ is \del or \delsave. %
		
		\smallskip
		
		\textbf{Case 1:} $\mathcal{O} = \del(G,f,u)$.
		
		\smallskip

		Then $G' = G - u$. Since $f_1(u) + f_2(u) = f(u) - 1 \geq -1$, we have $f_1(u) \geq 0$ or $f_2(u) \geq 0$. For concreteness, say $f_1(u) \geq 0$. Define a function $f_1' \colon V(G') \to \Z$ by
		\[
		f_1'(v) \,\defeq\, \begin{cases}
			f_1(v) - 1 &\text{if } uv \in E(G);\\
			f_1(v) &\text{otherwise}.
		\end{cases}
		\]
		Then $f_1' + f_2 = f' - 1$, so, by the inductive hypothesis, there is a partition $V(G') = V_1' \sqcup V_2$ such that $G[V_1']$ is weakly $f_1'$-degenerate and $G[V_2]$ is weakly $f_2$-degenerate. Set $V_1 \defeq V_1'\sqcup \set{u}$. We claim that the partition $V(G) = V_1 \sqcup V_2$ is as desired. Since $G[V_2]$ is weakly $f_2$-degenerate by assumption, we just need to argue that $G[V_1]$ is weakly $f_1$-degenerate. As $f_1(u) \geq 0$, the function $f_1$ is non-negative on $V_1$. Now we are done since $\del(G[V_1], f_1, u) = (G[V_1'],f_1')$ and $G[V_1']$ is weakly $f_1'$-degenerate.
		
		\smallskip
		
		\textbf{Case 2:} $\mathcal{O} = \delsave(G,f,u,w)$.
		
		\smallskip
		
		Again we have $G' = G - u$. It will be convenient to assume that $f_1(w)$, $f_2(w) \geq -1$. If this is not the case and, say, $f_1(w) < -1$, then we replace $f_1$ and $f_2$ by the functions $f_1^\ast$, $f_2^\ast \colon V(G) \to \Z$ given by $f_1^\ast(w) \defeq -1$, $f_2^\ast(w) \defeq f(w)$, and $f_i^\ast(v) \defeq f_i(v)$ for all $i \in \set{1,2}$ and $v \neq w$. %
		We can do this because every weakly $f_i^\ast$-degenerate subgraph of $G$ is also weakly $f_i$-degenerate. For $i = 2$ this follows from Lemma~\ref{lemma:monotone} since $f_2 \geq f_2^\ast$ by definition. %
		On the other hand, if a subgraph $H$ of $G$ is weakly $f_1^\ast$-degenerate, then $w \not\in V(H)$ since $f_1^\ast(w) < 0$. As $f_1^\ast$ and $f_1$ agree on all vertices except $w$, $H$ must be weakly $f_1$-degenerate as well. %

		Since this application of \delsave is legal, we have $f(u) > f(w)$, which implies that $f_1(u) > f_1(w)$ or $f_2(u) > f_2(w)$. For concreteness, say $f_1(u) > f_1(w)$. Define $f_1' \colon V(G') \to \Z$ by
		\[
		f_1'(v) \,\defeq\, \begin{cases}
			f_1(v) - 1 &\text{if } uv \in E(G) \text{ and } v \neq w;\\
			f_1(v) &\text{otherwise}.
		\end{cases}
		\]
		Then $f_1' + f_2 = f' - 1$, so, by the inductive hypothesis, there is a partition $V(G') = V_1' \sqcup V_2$ such that $G[V_1']$ is weakly $f_1'$-degenerate and $G[V_2]$ is weakly $f_2$-degenerate. Set $V_1 \defeq V_1'\sqcup \set{u}$. We claim that the partition $V(G) = V_1 \sqcup V_2$ is as desired. We just need to argue that $G[V_1]$ is weakly $f_1$-degenerate. Since $f_1(u) > f_1(w) \geq -1$, we have $f_1(u) \geq 0$. Hence, $f_1$ is non-negative on $V_1$. It remains to observe that by a legal application of one of the operations \del, \delsave it is possible to reduce the pair $(G[V_1], f_1)$ to $(G[V_1'], f_1')$. Indeed, if $w \not \in V_1$, then $(G[V_1'],f_1') = \del(G[V_1], f_1, u)$, while if $w \in V_1$, then $(G[V_1'],f_1') = \delsave(G[V_1], f_1, u, w)$.
	\end{scproof}
	
	\section{On-line DP-coloring and weak degeneracy}\label{sec:paint}
	
	As mentioned in the introduction, DP-paint number is an on-line version of DP-chromatic number introduced by Kim, Kostochka, Li, and Zhu in \cite{Paint}. It is defined by means of a certain game on a graph $G$:
	
	\begin{defn}[\textls{DP-painting game}]
		Let $G$ be a graph and let $g \colon V(G) \to \N$ be a function. The \emphd{DP-painting game on $(G,g)$} is played between two players---\Lister and \Painter---as follows. The game proceeds in \emphd{rounds}, starting with Round $0$. At the start of Round $i$, we have a graph $G_i$, where we initially set $G_0 \defeq G$. \Lister then picks a list $L_i(u)$ of colors for each vertex $u \in V(G_i)$ and assigns to every edge $uv \in E(G_i)$ a matching $C_{i,uv}$ from $L_i(u)$ to $L_i(v)$ \ep{the matching $C_{i,uv}$ need not be perfect and, in particular, may be empty}. In response, \Painter picks a function $\phi_i$ defined on some subset $U_i \subseteq V(G_i)$ with the following properties:
		\begin{itemize}
			\item $\phi_i(u) \in L_i(u)$ for all $u \in U_i$ \ep{in particular, $L_i(u) \neq \0$ for all $u \in U_i$}; and
			\item $\phi_i(u) \phi_i(v) \not \in C_{i, uv}$ for all $u$, $v \in U_i$ that are adjacent in $G_i$. 
		\end{itemize}
		Then we set $G_{i+1} \defeq G_i - U_i$ and proceed to Round $i+1$. \Lister wins the game if for some $i \in \N$, there is a vertex $u \in V(G_{i})$ with $\sum_{j < i} |L_j(u)| \geq g(u)$; otherwise, \Painter wins.

		\smallskip
		
		\noindent A graph $G$ is \emphd{$g$-DP-paintable} if \Painter has a winning strategy in the DP-painting game on $(G,g)$. Given $k \in \N$, we say that $G$ is \emphd{$k$-DP-paintable} if it is DP-paintable with respect to the constant $k$ function. The \emphd{DP-paint number} $\chi_{DPP}(G)$ of $G$ is the least $k$ such that $G$ is $k$-DP-paintable.
	\end{defn}
	
	Take $k \in \N$ and consider the DP-painting game on $(G, \bm{k})$, where $\bm{k}$ is the constant $k$ function. On Round $0$, \Lister may decide to give each vertex $u \in V(G)$ a list $L_0(u)$ of colors of size $|L_0(u)| = k$. Then \Painter must immediately assign a color to every vertex. Therefore, \Painter can win only if $\chi_{DP}(G) \leq k$, which shows that $\chi_{DP}(G) \leq \chi_{DPP}(G)$ for all $G$. On the other hand, if \Lister always plays so that $|L_i(u)| \leq 1$ for all $i$ and $u \in V(G_i)$, then \Painter can win if and only if $\chi_P(G) \leq k$, where $\chi_P(G)$ is the classical \emph{paint number} of $G$, i.e., the on-line analog of list-chromatic number \ep{see \cite[\S2]{Paint} for details}. Thus, $\chi_P(G) \leq \chi_{DPP}(G)$ as well, so $\chi_{DPP}(G)$ provides a common upper bound on $\chi_{DP}(G)$ and $\chi_P(G)$. It is shown in \cite{Paint} that either inequality $\chi_{DP}(G) \leq \chi_{DPP}(G)$ and $\chi_P(G) \leq \chi_{DPP}(G)$ can be strict; however, it is unknown if both of them can be strict at the same time. It is also not known if the difference $\chi_{DPP}(G) - \chi_{DP}(G)$ can be arbitrarily large.
	
	The goal of this section is to prove Proposition~\ref{prop:paint}, which says that the DP-paint number is bounded above by weak degeneracy plus $1$. We prove it in the following stronger form: %

	\begin{propcopy}{prop:paint}
		If $G$ is a weakly $f$-degenerate graph, then $G$ is $(f+1)$-DP-paintable.
	\end{propcopy}
	\begin{scproof}
		The strategy for \Painter is to pick functions $\phi_i$ so as to maintain the following property:
		
		\smallskip
		
		\begin{quote}
			\ep{$\mathrm{P}_i$}\quad $G_i$ is weakly $f_i$-degenerate, where $f_i(u) \defeq f(u) - \sum_{j < i} |L_j(u)|$ for all $u \in V(G_i)$. %
		\end{quote}
		
		\smallskip
		
		\noindent If this can be achieved, then \Painter will never lose, since for all $u \in V(G_i)$, we will have $f_i(u) \geq 0$, or, equivalently, $\sum_{j < i} |L_j(u)| < f(u) + 1$. Since $f_0 = f$, property \ep{$\mathrm{P}_0$} holds by assumption, so it remains to argue that if \ep{$\mathrm{P}_i$} holds at the start of Round $i$, then \Painter will be able to pick $\phi_i$ so that \ep{$\mathrm{P}_{i+1}$} holds.
		
		Suppose \Lister assigned a list $L_i(u)$ of colors and a matching $C_{i,uv}$ to each vertex $u \in V(G_i)$ and edge $uv \in E(G_i)$ respectively. For all $u \in V(G_i)$, let \[f_{i,1}(u) \,\defeq\, |L_i(u)| - 1 \quad \text{and} \quad f_{i,2}(u) \,\defeq\, f_i(u) - |L_i(u)|.\] Since $f_{i,1}(u) + f_{i,2}(u) = f_i(u) - 1$ for all $u \in V(G_i)$, Lemma~\ref{lemma:part} yields a partition $V(G_i) = U_i \sqcup W_i$ such that $G[U_i]$ is weakly $f_{i,1}$-degenerate, while $G[W_i]$ is weakly $f_{i,2}$-degenerate. By Lemma~\ref{lemma:DP}, $G[U_i]$ admits a proper $(L_i,C_i)$-coloring $\phi_i$. \Painter plays this coloring $\phi_i$. Then $G_{i+1} = G[W_i]$, so, to establish \ep{$\mathrm{P}_{i+1}$}, we need to show that $G[W_i]$ is weakly $f_{i+1}$-degenerate. To this end, note that for each $u \in W_i$,
		\[
		f_{i+1}(u) \,=\, f(u) - \sum_{j \leq i} |L_j(u)| \,=\, f_i(u) - |L_i(u)| \,=\, f_{i,2}(u),
		\]
		and $G[W_i]$ is indeed weakly $f_{i,2}$-degenerate by construction.
	\end{scproof}

	\section{Planar graphs are weakly 4-degenerate}\label{sec:planar}
	
	In this section we prove the analog of Thomassen's theorem \cite{Thom} on $5$-list-colorability of planar graphs in the context of weak degeneracy:
	
	\begin{theocopy}{theo:planar}\label{theo:planar1}
		Every planar graph is weakly 4-degenerate.
	\end{theocopy}
	
	As in the proof of Thomassen's theorem, we use induction to establish a technical lemma, which then easily yields Theorem~\ref{theo:planar1}. First, we need a definition. Let $G$ be a graph and let $f \colon V(G) \to \N$ be a function. Given a subset $U \subseteq V(G)$, we say that $G$ is \emphd{$U$-safely weakly $f$-degenerate} if, starting with $(G, f)$, it is possible to remove all vertices from $G$ by a sequence of legal applications of the operations \del and \delsave, so that every vertex in $U$ is removed using the \del operation. In particular, $G$ is $V(G)$-safely weakly $f$-degenerate if and only if $G$ is $f$-degenerate.

	\begin{lemma}\label{lemma:planar}
		Let $G$ be a planar graph on at least $3$ vertices where every non-outer face is triangular and the outer face is a cycle $C$ of length $k$. Let the vertices of $C$ in the natural order be $v_1$, \ldots, $v_k$. Define $f \colon V(G) \setminus \{v_1, v_2\} \to \mathbb{Z}$ by
		\[f(u)\,\defeq\, \begin{cases}
			2 - \abs{N_G(u) \cap \{v_1, v_2\}} & \text{if $u \in V(C)$}; \\
			4 - \abs{N_G(u) \cap \{v_1, v_2\}} & \text{otherwise}.
		\end{cases}\]
		Then $G - v_1- v_2$ is $(V(C) \setminus \set{v_1, v_2})$-safely weakly $f$-degenerate.
	\end{lemma}
	\begin{scproof}
		We proceed by induction on $\abs{V(G)}$. If $\abs{V(G)} = 3$, then $G - v_1 - v_2$ comprises a single vertex, which is $0$-degenerate, as desired. Now suppose that $\abs{V(G)} \geq 4$ and that the induction hypothesis holds for smaller graphs. We consider two cases.
		
		\smallskip
		
		\textbf{Case 1:} $C$ has a chord $v_a v_b$.
		
		\smallskip
		
		Then $C + v_a v_b$ is the union of two cycles $C_1$, $C_2$ with $E(C_1) \cap E(C_2) = \set{v_a v_b}$. Without loss of generality, suppose $v_1v_2 \in E(C_1)$ (and so $v_1v_2 \not\in E(C_2)$). Let $G_1$, $G_2$ be the respective induced subgraphs of $G$ on the vertices of each $C_i$ along with the vertices on the interiors of each cycle. Let $f_1 \defeq \rest{f}{V(G_1) \setminus \{v_1, v_2\}}$ and define $f_2 \colon V(G_2) \setminus \set{v_1, v_2} \to \Z$ by
		\[
		f_2(u) \,\defeq\, \begin{cases}
			2 - \abs{N_G(u) \cap \{v_a, v_b\}} & \text{if $u \in V(C_2)$}; \\
			4 - \abs{N_G(u) \cap \{v_a, v_b\}} & \text{otherwise}.
		\end{cases}
		\]
		By the induction hypothesis, starting with $(G_1 - v_1 - v_2, f_1)$, we can remove all vertices from $V(G_1) \setminus \set{v_1, v_2}$ via legal applications of the operations \del and \delsave, where each vertex in $V(C_1)\setminus \set{v_1, v_2}$ is removed using \del. Applying the same sequence of operations but starting with $(G - v_1 - v_2, f)$ yields the pair $(G_2 - v_a - v_b, f_2)$. By the inductive hypothesis again, we can now remove every remaining vertex via a sequence of legal applications of \del and \delsave, with every vertex in $V(C_2) \setminus \set{v_a, v_b}$ removed using \del, as desired.
		
		\smallskip
		
		\textbf{Case 2:} $C$ has no chord.
		
		\smallskip
		
		Since every non-outer face of $G$ is a triangle, the neighbors of $v_k$ form a path $u_1\ldots u_\ell$, where $u_1 = v_1$ and $u_\ell = v_{k-1}$. The assumption that $C$ has no chord implies that $u_2$, \ldots, $u_{\ell - 1}$ belong to the interior of $C$. Then the cycle $C' \defeq u_1\ldots u_{\ell} v_{k-2} \ldots v_1$ bounds the outer face of $G' \defeq G - v_k$. Applying the induction hypothesis to $G'$ shows that $G'' \defeq G - v_1 - v_2 - v_k$ is $(V(C')\setminus \set{v_1, v_2})$-safely weakly $f'$-degenerate, where $f' \colon V(G'') \to \mathbb{Z}$ is defined by
		\[f'(u)\,\defeq\, \begin{cases}
			2 - \abs{N_{G'}(u) \cap \{v_1, v_2\}} & \text{if $u \in V(C')$}; \\
			4 - \abs{N_{G'}(u) \cap \{v_1, v_2\}} & \text{otherwise}.
		\end{cases}.\]
		In other words, starting with $(G'', f')$, we can remove every vertex by a sequence of legal applications of \del and \delsave, where each vertex in $V(C')\setminus \set{v_1, v_2}$ is removed using \del. %
		Since $f' \leq f$, we may apply the same sequence of operations starting with $(G - v_1 - v_2, f)$ instead \ep{see Lemma~\ref{lemma:monotone}}. Moreover, we can accrue some extra savings for the vertex $v_k$, as follows. Consider any $u_i$ with $2 \leq i \leq \ell - 1$. By assumption, $u_i$ is removed from $G''$ using the \del operation, but since $u_i \not \in V(C)$, we are now allowed to remove $u_i$ using \delsave. Notice that $f'(u_i) = f(u_i) - 2$, because $u_i$ is in $V(C')$ but not in $V(C)$. When $u_i$ was removed from $G''$, the value of the function at $u_i$ was at least $0$, which means that at the same stage of the process on $G - v_1 - v_2$, the value of the function at $u_i$ is at least $2$. On the other hand, since $v_k \in V(C)$ and $v_1 \in N_G(v_k)$, we have $f(v_k) \leq 1 < 2$. This means that instead of using the operation \del, we can legally remove $u_i$ using $\delsave(\cdot,\cdot, u_i, v_k)$. Upon performing this modified sequence of operations, we only have $v_k$ left to remove, so we just need to check that the value of the function at $v_k$ is at least $0$. To this end, note that $f(v_k)$ is $1$ if $v_{k-1} \neq v_2$ and $0$ otherwise. Since the only neighbor of $v_k$ that may be removed without saving $v_k$ is $v_{k-1}$, and that can only happen when $v_{k-1} \neq v_2$, it follows that the value at $v_k$ cannot drop below $0$, as desired.
	\end{scproof}
	
	We now complete the proof of the theorem.
	
	\begin{scproof}[ of Theorem~\ref{theo:planar1}]
		Since adding vertices or edges cannot decrease the weak degeneracy of a graph, it suffices to prove the theorem for maximal planar graphs $G$ on at least $3$ vertices. Then $G$ is a planar triangulation. Let $v_1$, $v_2$ be adjacent vertices on its outer face. Removing $v_1$ and $v_2$ using \del and then applying Lemma~\ref{lemma:planar}, we see that $G$ is weakly $4$-degenerate, as desired.
	\end{scproof}
	
	\section{Brooks-type results}\label{sec:GDP}
	
	\subsection{Weakly $(\deg-1)$-degenerate graphs}
	
	We say that a graph $G$ is \emphd{weakly $(\deg - 1)$-degenerate} if it is weakly degenerate with respect to the function $u \mapsto \deg_G(u) - 1$. Recall that a \emphd{GDP tree} is a connected graph in which every block is either a clique or a cycle.
	The main result of this section is the following characterization of connected weakly $(\deg - 1)$-degenerate graphs:
	
	\begin{theocopy}{theo:GDP}
		Let $G$ be a connected graph. The following statements are equivalent:
		\begin{enumerate}[label=\ep{\normalfont\arabic*}]
			\item\label{item:deg-1} $G$ is weakly $(\deg - 1)$-degenerate;
			\item\label{item:tree} $G$ is not a GDP-tree.
		\end{enumerate}
	\end{theocopy}
	
	To begin with, we need the following standard fact:
	
	\begin{lemma}\label{lemma:more}
		Let $G$ be a connected graph and let $f \colon V(G) \to \N$ be a function. Suppose that:
		\begin{enumerate}[label=\ep{\itshape\alph*}]
			\item\label{item:a} $f(u) \geq \deg_G(u) - 1$ for all $u \in V(G)$; and
			\item\label{item:b} $f(x) \geq \deg_G(x)$ for some $x \in V(G)$.
		\end{enumerate}
		Then $G$ is $f$-degenerate.
	\end{lemma}
	\begin{scproof}
		Fix a vertex $x$ witnessing \ref{item:b} and list the vertices of $G$ as $u_1$, $u_2$, \ldots, $u_n$ in order of decreasing distance to $x$, resolving ties arbitrarily. Then $u_n = x$ and, for each $1 \leq i \leq n-1$, the vertex $u_i$ has at least one neighbor among $u_{i+1}$, \ldots, $u_n$. We can now remove all vertices from $G$ by applying the operation \del to them in this order.
	\end{scproof}
	
	The next lemma contains the central part of our argument:
	
	\begin{lemma}\label{lemma:2con}
		Let $G$ be a connected graph that is not weakly $(\deg-1)$-degenerate. Then every connected induced subgraph of $G$ without cut vertices is regular.
	\end{lemma}
	\begin{scproof}
		Take a subset $A \subseteq V(G)$ such that the subgraph $G[A]$ has no cut vertices and suppose, toward a contradiction, that $G[A]$ is not regular. Define $f \colon V(G) \to \N$ by $f(u) \defeq \deg_G(u) - 1$ for all $u \in V(G)$. Our goal is to show that $G$ is weakly $f$-degenerate. Note that every connected component of $G - A$ contains at least one vertex $v$ that has a neighbor in $A$ and hence satisfies $f(v) \geq \deg_{G - A}(v)$. Therefore, by Lemma~\ref{lemma:more}, we can remove all vertices from $G - A$ using only the operation \del. After this, the graph $G$ will be replaced by $G' \defeq G[A]$ and the function $f$ by $f' \colon A \to \N \colon u \mapsto \deg_{G[A]}(u) - 1$. Since the graph $G[A]$ is connected and not regular, we can pick two adjacent vertices $x$, $y \in A$ with $\deg_{G[A]}(x) < \deg_{G[A]}(y)$ and hence $f'(x) < f'(y)$. Now we let
		\[
		(G'', f'') \,\defeq\, \delsave(G[A],f', y, x).
		\]
		Since $f'(y) > f'(x)$, this is a legal application of $\delsave$. As the graph $G[A]$ has no cut vertices, the graph $G'' = G[A] - y$ is connected. It remains to observe that $G''$ is $f''$-degenerate by Lemma~\ref{lemma:more}, where condition \ref{item:b} is witnessed by the vertex $x$.
	\end{scproof}
	
	It remains to characterize the graphs satisfying the conclusion of Lemma~\ref{lemma:2con}:
	
	\begin{lemma}\label{lemma:tree}
		Let $G$ be a connected graph such that every connected induced subgraph of $G$ without cut vertices is regular. Then $G$ is a GDP-tree.
	\end{lemma}
	\begin{scproof}
		Suppose, toward a contradiction, that $G$ is a counterexample with the fewest vertices. Note that $|V(G)| \geq 4$, since all connected graphs on at most $3$ vertices are GDP-trees. By the minimality of $|V(G)|$, every proper connected induced subgraph of $G$ must be a GDP-tree.
		
		We claim that $G$ is $2$-connected. Otherwise, every block in $G$ would be a proper connected induced subgraph of $G$, hence a GDP-tree. The only GDP-trees without cut vertices are cliques and cycles, so this implies that every block in $G$ is a clique or a cycle, i.e., $G$ is a GDP-tree.
		
		Since $G$ is $2$-connected, it must be regular. Let $d$ be the common degree of the vertices of $G$. Then $d \geq 2$ by $2$-connectedness. Furthermore, if $d$ were equal to $2$, then $G$ would be a cycle and hence a GDP-tree. Therefore, $d \geq 3$.
		
		Pick an arbitrary vertex $x \in V(G)$ and consider the graph $G' \defeq G - x$. Then $G'$ is connected, so it is a GDP-tree. Since $G$ is regular and not a clique, not every vertex in $V(G')$ is adjacent to $x$. This implies that $G'$ is not regular, so it must have a cut vertex and at least two blocks.
		
		Let $B$ be an arbitrary leaf block in $G'$ and let $c \in V(B)$ be the cut vertex of $G'$ in $B$. The graph $B$ is regular, so let $k$ be the common degree of every vertex of $B$. %
		The degree of a vertex $u \in V(B) \setminus \set{c}$ in $G$ is either $k + 1$ or $k$, depending on whether $u$ is adjacent to $x$ or not. Since $G$ is $2$-connected, $x$ must be adjacent to at least one vertex in $V(B) \setminus \set{c}$, which, since $G$ is $d$-regular, implies that $k+1 = d$ and $x$ is in fact adjacent to every vertex in $V(B) \setminus \set{c}$. Hence, $x$ has at least $|V(B)| - 1 \geq k = d - 1$ neighbors in $B$.
		
		Finally, as there are at least $2$ distinct leaf blocks in $G'$, we conclude that $x$ has at least $2(d-1)$ neighbors. Therefore, $d \geq 2(d-2)$, i.e., $d \leq 2$, which is a contradiction.
	\end{scproof}
	
	Theorem~\ref{theo:GDP} now follows easily:
	
	\begin{scproof}[ of Theorem~\ref{theo:GDP}]
		The implication \ref{item:tree} $\Longrightarrow$ \ref{item:deg-1} is a combination of Lemmas~\ref{lemma:2con} and \ref{lemma:tree}. The implication \ref{item:deg-1} $\Longrightarrow$ \ref{item:tree} follows since GDP-trees are not DP-degree-colorable \cite[Theorem 9]{BKP}. That is, if $G$ is a GDP-tree, then it is possible to give each vertex $u \in V(G)$ a list $L(u)$ of available colors of size $|L(u)| \geq \deg_G(u)$ and assign to each edge $uv \in E(G)$ a matching $C_{uv}$ from $L(u)$ to $L(v)$ so that $G$ is does not admit a proper $(L,C)$-coloring. By Lemma~\ref{lemma:DP}, this implies that $G$ is not weakly $(\deg - 1)$-degenerate.
	\end{scproof}
	
	\subsection{Weak degeneracy and maximum average degree}\label{sec:mad}
	
	Here we establish a lower bound on the maximum average degree of a graph in terms of its weak degeneracy:
	
	\begin{theocopy}{theo:mad}
		Let $G$ be a nonempty graph. If the weak degeneracy of $G$ is at least $d \geq 3$, then either $G$ contains a $(d+1)$-clique or
		\[
		\mathrm{mad}(G) \,\geq\, d \,+\, \frac{d-2}{d^2 + 2d - 2}.
		\]
	\end{theocopy}
	
	We derive Theorem~\ref{theo:mad} from Theorem~\ref{theo:GDP}. Our argument is closely analogous to the proof of the lower bound on the average degree of DP-critical graphs due to Kostochka, Pron, and the first named author \cite[Corollary 10]{BKP}, which in turn is based on earlier work of Gallai \cite{Gallai}.
	
	We need the following result, essentially established by Gallai in \cite{Gallai} \ep{Gallai's paper is in German; see \cite[Appendix]{BKP} for a proof in English}:
	
	\begin{lemma}[{\cite[Lemma 20]{BKP}}]\label{lemma:Gallai}
		Let $T$ be a GDP-tree of maximum degree at most $d \geq 3$ and without a $(d+1)$-clique. Then $\ad(T)\leq d - 1 + 2/d$.
	\end{lemma}
	
	We say that $G$ is a \emphd{minimal} graph of weak degeneracy $d$ if $\wdg(G) = d$ and $\wdg(H) < d$ for every proper subgraph $H$ of $G$.
	
	\begin{lemma}\label{lemma:crit}
		Let $G$ be a minimal graph of weak degeneracy $d \geq 3$. %
		\begin{enumerate}[label=\ep{\itshape\alph*}]
			\item\label{item:crita} The minimum degree of $G$ is at least $d$.
			
			\item\label{item:critb} Let $U \defeq \set{u \in V(G) \,:\, \deg_G(u) = d}$. Then every component of $G[U]$ is a GDP-tree.
		\end{enumerate}
	\end{lemma}
	\begin{scproof}
		\ref{item:crita} Suppose that there is a vertex $u \in V(G)$ with $\deg_G(u) \leq d - 1$. We will show that $G$ is weakly $(d-1)$-degenerate. Let $f$ be the constant $d-1$ function on $V(G)$. By the minimality of $G$, we may remove every vertex from $(G - u, f)$ via a sequence of legal applications of the operations \del and \delsave. Since $\deg_G(u) \leq d - 1$, we may use the same sequence of operations to remove every vertex except $u$ from $(G, f)$ \ep{at which point the function $f$ will be replaced by the map sending $u$ to $d - 1 - \deg_G(u)$} and then remove $u$ using the operation \del.
		
		\ref{item:critb} Let $C$ be a connected component of $G[U]$ and let $f$ be the constant $d - 1$ function on $V(G)$. By the minimality of $G$, we may remove every vertex from $(G - V(C), f)$ via a sequence of legal applications of the operations \del and \delsave. If we perform the same sequence of operations on $(G,f)$, then the graph $G$ will be replaced by $C$, while the function $f$ will be replaced by the map sending each $u \in V(C)$ to $d-1 - \deg_{G - V[C]}(u) = \deg_C(u) - 1$. Since $G$ is not weakly $(d-1)$-degenerate, this implies that $C$ is not weakly $(\deg - 1)$-degenerate. Hence, by Theorem~\ref{theo:GDP}, $C$ is a GDP-tree.
	\end{scproof}
	
	\begin{scproof}[ of Theorem~\ref{theo:mad}]
		Fix $d \geq 3$. It suffices to argue that every minimal graph $G$ of weak degeneracy $d$ and without a $(d+1)$-clique satisfies
		\[
		\ad(G) \,\geq\, d \,+\, \frac{d-2}{d^2 + 2d - 2}.
		\]
		To this end, we use discharging. Let the initial charge of each vertex $u \in V(G)$ be $\mathrm{ch}(u) \defeq \deg_G(u)$. The only discharging rule is: Every vertex $u \in V(G)$ with $\deg_G(u) > d$ sends to each neighbor the charge $c \defeq d/(d^2 + 2d - 2)$. Let the new charge of each vertex $u$ be $\mathrm{ch}^\ast(u)$. Note that
		\[
		\ad(G) |V(G)| \,=\, \sum_{u \in V(G)} \mathrm{ch}(u) \,=\, \sum_{u \in V(G)} \mathrm{ch}^\ast(u).
		\]
		For any vertex $u$ with $\deg_G(u) > d$, we have
		\[
		\mathrm{ch}^\ast(u) \,\geq\, \deg_G(u) - c \deg_G(u) \,\geq\, \left(1 \,-\, \frac{d}{d^2 + 2d -2}\right)(d+1) \,=\, d \,+\, \frac{d-2}{d^2 + 2d - 2}.
		\]
		Let $C$ be any connected component of $G[U]$, where $U$ is the set of all vertices of degree $d$ in $G$. By Lemma~\ref{lemma:crit}\ref{item:critb}, $C$ is a GDP-tree. Hence, by Lemma~\ref{lemma:Gallai}, $\ad(C)\leq  d -1 + 2/d$. Therefore, %
		\begin{align*}
			\sum_{u \in V(C)} \mathrm{ch}^\ast(u) \,&=\, d|V(C)| \,+\, c \sum_{u \in V(C)} (d - \deg_C(u)) \\
			&\geq\, \left(d \,+\, \frac{d}{d^2 + 2d -2} \left(1 - \frac{2}{d}\right)\right)|V(C)|\,=\, \left(d \,+\, \frac{d-2}{d^2 + 2d - 2}\right)|V(C)|.
		\end{align*}
		The above bounds imply that
		\[
		\sum_{u \in V(G)} \mathrm{ch}^\ast(u) \,\geq\, \left(d \,+\, \frac{d-2}{d^2 + 2d - 2}\right)|V(G)|,
		\]
		which yields the desired result.
	\end{scproof}
	
	\section{Lower bounds for regular graphs}\label{sec:Kdd}
	
	In this section we establish lower bounds on weak degeneracy for regular graphs.
	
	\begin{propcopy}{prop:lower_bound}\label{prop:lower_bound_copy}
		Let $G$ be a $d$-regular graph with $n \geq 2$ vertices. Then $\wdg(G) \geq d - \sqrt{2n}$.
	\end{propcopy}
	\begin{scproof}
		Let $k \defeq \wdg(G)$. Set $G_0 \defeq G$ and let $f_0$ be the constant $k$ function on $V(G)$. By definition, starting with $(G_0,f_0)$, it is possible to remove all vertices from $G$ via a sequence of legal applications of the operations \del and \delsave. Fix any such sequence $\mathcal{S} = (\mathcal{O}_0, \ldots, \mathcal{O}_{n-1})$. For each $0 \leq i \leq n-1$, let $(G_{i+1},f_{i+1})$ be the result of applying $\mathcal{O}_i$ to $(G_i, f_i)$. Then we can write
		\[
		\mathcal{O}_i \,=\, \del(G_i, f_i, u_i) \quad \text{or} \quad \mathcal{O}_i \,=\, \delsave(G_i,f_i,u_i,w_i).
		\]
		For each $0 \leq i \leq n-1$, define
		\[
		d_i\,\defeq\, |\set{j < i \,:\, u_ju_i \in E(G)}| \quad \text{and} \quad \sigma_i \,\defeq\, |\set{j < i \,:\, \mathcal{O}_j = \delsave(G_j,f_j,u_j,u_i)}|.
		\]
		\ep{So $\sigma_i$ is the number of vertices that ``save'' $u_i$.} Then $0 \leq f_i(u_i) = k - d_i + \sigma_i$ and thus
		\begin{equation}\label{eq:k}
			k \,\geq\, d_i - \sigma_i.
		\end{equation}
		Adding \eqref{eq:k} up over the interval $n-t \leq i \leq n-1$ for some integer $1 \leq t \leq n$ yields
		\begin{equation}\label{eq:kt}
			kt \,\geq\, \left(\sum_{i = n - t}^{n-1} d_i\right) \,-\, \left(\sum_{i = n - t}^{n-1} \sigma_i\right).
		\end{equation}
		Each index $j$ contributes to $\sigma_i$ for at most one $i$, so $\sum_{i = 0}^{n-1} \sigma_i \leq n$. Also, since $G$ is $d$-regular,
		\begin{equation}\label{eq:deg_sum}
			\sum_{i = n - t}^{n-1} d_i \,=\, dt - |E(G[u_{n-t}, \ldots, u_{n-1}])| \,\geq\, dt - {t \choose 2}.
		\end{equation}
		Therefore, \eqref{eq:kt} implies that
		\[
		k \,\geq\, d - \frac{t-1}{2} - \frac{n}{t}.
		\]
		Finally, taking $t \defeq \lceil \sqrt{2n} \rceil$ gives
		\[
		k \,\geq\, d - \frac{\lceil \sqrt{2n} \rceil - 1}{2} - \frac{n}{\lceil \sqrt{2n} \rceil} \,\geq\, d - \frac{\sqrt{2n}}{2} - \frac{n}{\sqrt{2n}} \,=\, d - \sqrt{2n},
		\]
		as desired.
	\end{scproof}
	
	\begin{propcopy}{prop:Kdd}
		If $G$ is a triangle-free $d$-regular graph with $n \geq 4$ vertices, then $\wdg(G) > d - \sqrt{n} - 1$.
	\end{propcopy}
	\begin{scproof}
		The argument is virtually the same as in the proof of Proposition~\ref{prop:lower_bound_copy}, except that we use Mantel's theorem to replace the bound \eqref{eq:deg_sum} by $\sum_{i = n - t}^{n-1} d_i \geq dt - t^2/4$. This yields
		\[
		\wdg(G) \,\geq\, d - \frac{t}{4} - \frac{n}{t}.
		\]
		Now we set $t \defeq \lceil 2\sqrt{n}\rceil$ and conclude that
		\[
		\wdg(G) \,\geq\, d - \frac{\lceil 2\sqrt{n}\rceil}{4} - \frac{n}{\lceil 2\sqrt{n}\rceil} \,\geq\, d - \frac{2\sqrt{n} + 1}{4} - \frac{n}{2\sqrt{n}} \,\geq\, d - \sqrt{n} -\frac{1}{4},
		\]
		as desired.
	\end{scproof}
	
	\section{Going below the maximum degree}\label{sec:chrom}
	
	\subsection{Preliminaries}
	
	In this section we review some necessary background facts. First, we will need the Lov\'asz Local Lemma, in the following form:
	
	\begin{theo}[{\textls{Lov\'asz Local Lemma} \cite[Lemma 5.1.1]{prob_method}}]\label{theo:LLL}
		Let $\mathcal{X}$ be a finite family of random events such that each $X \in {\mathcal{X}}$ has probability at most $p$ and is mutually independent from all but $\Delta$ other events in $\mathcal{X}$. If $e p (\Delta + 1) \leq 1$, then the probability that no event in $\mathcal{X}$ happens is positive.
	\end{theo}
	
	We shall also use the Chernoff bound for binomial random variables:
	
	\begin{theo}[{\textls{Chernoff bound} \cite[43]{MolloyReed}}]\label{theo:Chernoff}
		If $X \sim \mathrm{Bin}(n,p)$ is a binomial random variable, then for all $0 \leq t \leq np$,
		\[
		\P[|X - np| > t] \,<\, 2\exp\left(-\frac{t^2}{3np}\right).
		\]
	\end{theo}
	
	Next, we need a quantitative version of the Central Limit Theorem due to Berry and Esseen:
	
	\begin{theo}[{Berry--Esseen \cite[\S{}XVI.5]{Feller}}]\label{theo:BE}
		There is a universal constant $A > 0$ with the following property. Let $X_1$, \ldots, $X_n$ be independent identically distributed random variables such that $\E[X_i] = 0$, $\E[X_i^2] = \sigma^2 > 0$, and $\E[|X_i|^3] = \rho < \infty$. Then for all $t \in \R$,
		\[
		\left|\P\left[\frac{\sum_{i=1}^n X_i}{\sigma\sqrt{n}} \,\leq \, t \right] \,-\, \frac{1}{\sqrt{2\pi}}\int_{-\infty}^t e^{-x^2/2} \diff x\right| \,\leq \, \frac{A\rho}{\sigma^3 \sqrt{n}}.
		\]
	\end{theo}
	
	In particular, if $X \sim \mathrm{Bin}(n,p)$ is a binomial random variable, then for any $\beta > 0$,
	\begin{equation}\label{eq:BE}
		\P[X \leq np - \beta \sqrt{n}] \,=\, \int_{-\infty}^{\beta/\sqrt{p(1-p)}} e^{-x^2/2} \diff x \,+\, O\left(\frac{1-2p(1-p)}{\sqrt{p(1-p)n}}\right)
	\end{equation}
	This means that for large $n$, $\P[X \leq np-\beta\sqrt{n}]$ is separated from $0$. By applying this result to the random variable $n - X \sim \mathrm{Bin}(n, 1-p)$, we see that $\P[X \geq np + \beta\sqrt{n}]$ is separated from $0$ as well.
	
	The following is a standard consequence of Hall's theorem:
	
	\begin{lemma}\label{lemma:Hall}
		Let $G$ be a graph and let $A$, $B \subseteq V(G)$ be disjoint sets. Suppose that each vertex in $A$ has at most $d_1$ neighbors in $B$, while each vertex in $B$ has at least $d_2$ neighbors in $A$. Let $t \in \N$ satisfy $t d_1 \leq d_2$. Then there exists a partial function $s \colon A \pto B$ such that:
		\begin{itemize}
			\item for all $u \in A$, if $s(u)$ is defined, then $s(u)$ is a neighbor of $u$;
			\item the preimage of every vertex $w \in B$ under $s$ has cardinality exactly $t$.
		\end{itemize}
	\end{lemma}
	\begin{scproof}
		Let $H$ be the maximal bipartite subgraph of $G$ with parts $A$ and $B$, and let $H^\ast$ be obtained from $H$ by replacing every vertex $w \in B$ by $t$ copies, denoted $w_1$, \ldots, $w_t$. By construction, $H^\ast$ is a bipartite graph with parts $A$ and $B^\ast \defeq \set{w_j \,:\, w \in B, \, 1 \leq j \leq t}$. For all $u \in A$, $\deg_{H^\ast}(u) \leq t  d_1 \leq d_2$. On the other hand, every vertex $w_j \in B^\ast$ satisfies $\deg_{H^\ast}(w_j) \geq d_2$. These inequalities, together with Hall's theorem \cite[Theorem 2.1.2]{Die}, imply that $H^\ast$ has a matching $M$ that saturates $B^\ast$. We can now define the desired function $s \colon A \pto B$ by mapping each $u \in A$ that is covered by $M$ to the unique $w \in B$ such that $uw_j \in M$ for some $j$. 
	\end{scproof}
	
	It will be convenient for us to work with $d$-regular graphs rather than graphs of maximum degree $d$. To this end, we shall employ the following facts:
	
	\begin{lemma}[{Chartrand--Wall \cite{embedding}}]\label{lemma:embchrom}
		If $G$ is a graph of maximum degree $d$ and chromatic number at most $k$, then $G$ can be embedded into a $d$-regular graph $G^\ast$ of chromatic number at most $k$.
	\end{lemma}
	
	\begin{lemma}[{\cite[Exercise 12.4]{MolloyReed}}]\label{lemma:embgirth}
		If $G$ is a graph of maximum degree $d$ and girth at least $g$, then $G$ can be embedded into a $d$-regular graph $G^\ast$ of girth at least $g$.
	\end{lemma}
	\begin{scproof}
		This fact is well-known, but we include a proof for completeness. We use a simplified version of the construction from \cite[Proposition 4.1]{BipartiteNibble}. Set
		\[
		N \,\defeq\, \sum_{u \in V(G)}(d-\deg_G(u))
		\]
		and let $\Gamma$ be an $N$-regular graph of girth at least $g$, which exists by \cite{GirthRegular1, GirthRegular2}. We may assume that $V(\Gamma) = \set{1,\ldots,q}$, where $q \defeq |V(\Gamma)|$. Take $q$ vertex-disjoint copies of $G$, say $G_1$, \ldots, $G_q$ and define $S_i \defeq \set{u\in V(G_i)\,:\, \deg_{G_i}(u) < d}$ for every $1\leq i \leq q$. The graph $G^\ast$ is obtained from the disjoint union of $G_1$, \ldots, $G_q$ by performing the following sequence of operations once for each edge $ij \in E(\Gamma)$, one edge at a time:
		\begin{enumerate}[label=\ep{\arabic*}]
			\item Pick arbitrary vertices $u \in S_i$ and $v \in S_j$.
			\item Add the edge $uv$ to $E(G^\ast)$.
			\item If $\deg_{G^*}(u) = d$, remove $u$ from $S_i$.
			\item If $\deg_{G^*}(v) = d$, remove $v$ from $S_j$.
		\end{enumerate}
		It is clear that the resulting graph $G^\ast$ is as desired.
	\end{scproof}
	
	\subsection{Removal schemes}
	
	In the next definition we introduce the technical notion of a \emph{removal scheme} on a graph $G$. Roughly speaking, a removal scheme records the order in which we attempt to remove the vertices from $G$. Additionally, it indicates whether each vertex is removed using a \del or a \delsave operation, and in the latter case, what other vertex we ``save'' an extra color for.
	
	\begin{defn}[\textls{Removal schemes}]
		Fix a graph $G$. A \emphd{removal scheme} on $G$ is a pair $(\prec, \mathsf{save})$, where $\prec$ is a linear ordering on $V(G)$ and $\mathsf{save} \colon V(G) \pto V(G)$ is a partial function such that for every vertex $u \in V(G)$, if $\mathsf{save}(u)$ is defined, then it is a neighbor of $u$ and $u \prec \mathsf{save}(u)$. For convenience, we write $\mathsf{save}(u) = \mathsf{blank}$ if $\mathsf{save}(u)$ is undefined. Given a removal scheme $(\prec, \mathsf{save})$, we call $\prec$ the \emphd{removal order} and say that a vertex $u$ \emphd{saves} the vertex $\mathsf{save}(u)$. A removal scheme $(\prec, \mathsf{save})$ is \emphd{legal} if for all $u$, $w \in V(G)$ such that $w = \mathsf{save}(u)$, we have
		\begin{equation}\label{eq:legal}
			|\set{v \in N_G(u) \,:\, v \prec u \text{ and } \mathsf{save}(v) \neq u}| \,<\, |\set{v \in N_G(w) \,:\, v \prec u \text{ and } \mathsf{save}(v) \neq w}|.
		\end{equation}
		The \emphd{gap} of a vertex $u$ with respect to a removal scheme $(\prec, \mathsf{save})$ is the quantity
		\[
		\mathsf{gap}(u; \prec, \mathsf{save}) \,\defeq\, |\set{v \in N_G(u) \,:\, v \succ u}| \,+\, |\set{v \in N_G(u) \,:\, \mathsf{save}(v) = u}|.
		\]
		We also let $\mathsf{gap}(\prec, \mathsf{save}) \defeq \min \set{\mathsf{gap}(u; \prec, \mathsf{save}) \,:\, u \in V(G)}$.
	\end{defn}
	
	\begin{lemma}\label{lemma:remove}
		Let $G$ be a graph of maximum degree at most $d$ and let $(\prec, \mathsf{save})$ be a legal removal scheme on $G$. Then $G$ is weakly $(d - \mathsf{gap}(\prec, \mathsf{save}))$-degenerate.
	\end{lemma}
	\begin{scproof}
		For brevity, let $g \defeq \mathsf{gap}(\prec, \mathsf{save})$. Let $u_0$, \ldots, $u_{n-1}$ be the vertices of $G$ listed in the order given by $\prec$. Define a sequence $(G_i, f_i)$, $0 \leq i \leq n-1$ by setting $(G_0, f_0) \defeq (G, d - g)$ and
		\[
		(G_{i+1},f_{i+1}) \,\defeq\, \begin{cases}
			\del(G_i, f_i, u_i) &\text{if } \mathsf{save}(u_i) = \mathsf{blank};\\
			\delsave(G_i,f_i, u_i, \mathsf{save}(u_i)) &\text{otherwise}.
		\end{cases}
		\]
		We claim that this construction yields a sequence of legal applications of \del and \delsave that removes every vertex from $G$. Indeed, consider any vertex $u_i$. Note that
		\[
		f_i(u_i) \,=\, d - g - |\set{v \in N_G(u_i) \,:\, v \prec u_i \text{ and } \mathsf{save}(v) \neq u_i}| \,\geq\, \mathsf{gap}(u_i;\prec\mathsf{save}) - g \,\geq\, 0. 
		\]
		This shows that the functions $f_i$ are non-negative. Now suppose that $\mathsf{save}(u_i) = w \in V(G)$. Then, by definition, $w$ is a neighbor of $u_i$ that appears after $u_i$ in the ordering $\prec$, and thus the operation $\delsave(G_i,f_i, u_i, w)$ may be applied. Furthermore, $f_i(u_i) > f_i(w)$ by \eqref{eq:legal}, so this application of $\delsave$ is legal, as desired.
	\end{scproof}
	
	\subsection{Regular sets}
	
	Let $G$ be a graph. Given a vertex $u \in V(G)$ and a set $A \subseteq V(G)$, we let $N_A(u) \defeq N_G(u) \cap A$ denote the set of all neighbors of $u$ in $A$ and write $\deg_A(u) \defeq |N_A(u)|$. Several times in our arguments, we will need to perform the following operation: given a set $A$ and a number $p \in [0,1]$, we will need to pick a subset $A' \subseteq A$ such that every vertex $u \in V(G)$ has roughly $p\deg_A(u)$ neighbors in $A'$. Formally, we introduce the following definition:
	
	\begin{defn}[\textls{Regular sets}]
		Fix a graph $G$ of maximum degree $d$ and a subset $A \subseteq V(G)$. Given $p$, $\epsilon \in (0,1]$, a \emphd{$(p,\epsilon)$-regular subset of $A$} is a set $A' \subseteq A$ such that every vertex $u \in V(G)$ satisfies one of the following conditions:
		\begin{itemize}
			\item either $\deg_A(u) < 9\log d/(\epsilon^2 p)$ \ep{i.e., $u$ has very few neighbors in $A$},
			\item or $(1-\epsilon) p\deg_A(u) \leq \deg_{A'}(u) \leq (1+\epsilon)p\deg_A(u)$ \ep{i.e., $\deg_{A'}(u) \approx p\deg_A(u)$}.
		\end{itemize}
	\end{defn}
	
	Using the Lov\'asz Local Lemma, it is not hard to prove that $(p,\epsilon)$-regular subsets exist:
	
	\begin{lemma}\label{lemma:regular}
		Let $G$ be a graph of maximum degree $d$. Fix $p$, $\epsilon \in (0,1]$. Then every set $A \subseteq V(G)$ has a $(p,\epsilon)$-regular subset.
	\end{lemma}
	\begin{scproof}
		We may assume $d > 8$, as otherwise $\deg_A(u) \leq d < 9\log d$ for all $u \in V(G)$, so any subset $A' \subseteq A$ is $(p, \epsilon)$-regular. Form a random set $A' \subseteq A$ by picking each vertex independently with probability $p$. We shall use the \hyperref[theo:LLL]{Lov\'asz Local Lemma} \ep{Theorem~\ref{theo:LLL}} to argue that $A'$ is $(p,\epsilon)$-regular with positive probability. Let $U \subseteq V(G)$ be the set of all vertices $u \in V(G)$ with $\deg_A(u) \geq 9\log d/(\epsilon^2 p)$. For each $u \in U$, let $X_u$ be the random event that
		\[
		\deg_{A'}(u) \,\not\in \left[(1-\epsilon) p\deg_A(u),\, (1+\epsilon)p\deg_A(u)\right].
		\]
		We need to argue that with positive probability, none of the events $X_u$ happen. By the \hyperref[theo:Chernoff]{Chernoff bound} \ep{Theorem~\ref{theo:Chernoff}}, for each $u \in U$ we have
		\[
		\P[X_u] \,<\, 2 \exp\left(-\frac{\epsilon^2 p \deg_A(u)}{3}\right) \,\leq\, 2\exp\left(-3\log d\right) \,=\, 3d^{-3}. 
		\]
		Each event $X_u$ is mutually independent from all the events $X_v$ corresponding to the vertices $v$ that do not share a neighbor with $u$. Since there are at most $d(d-1)$ vertices that share a neighbor with $u$ \ep{not including $u$ itself}, the \hyperref[theo:LLL]{Lov\'asz Local Lemma} shows that with positive probability none of the events $X_u$, $u \in U$ happen provided that
		\[
		e \cdot 3d^{-3} \cdot (d(d-1) + 1) \,<\, 1.
		\]
		This inequality holds for all $d > 8$, and the proof is complete.
	\end{scproof}
	
	\subsection{Graphs of bounded chromatic number}

	\begin{theocopy}{theo:chrom}\label{theo:chrom1}
		For each integer $k \geq 1$, there exist $c > 0$ and $d_0 \in \N$ such that if $G$ is a graph of maximum degree $d \geq d_0$ with $\chi(G) \leq k$, then $\mathsf{wd}(G) \leq d - c\sqrt{d}$.
	\end{theocopy}
	
	Let $G$ be a graph of maximum degree $d$ and chromatic number at most $k$, where we assume that $d$ is sufficiently large in terms of $k$. Upon replacing $G$ with a supergraph if necessary, we may assume that $G$ is $d$-regular \ep{Lemma~\ref{lemma:embchrom}}. Let $c$ be a sufficiently small positive quantity depending on $k$ \ep{but not on $d$}. We will construct a legal removal scheme $(\prec, \mathsf{save})$ on $G$ such that $\mathsf{gap}(\prec, \mathsf{save}) \geq c\sqrt{d}$. By Lemma~\ref{lemma:remove}, this will yield the desired result.
	
	We start by applying Lemma~\ref{lemma:regular} to obtain a $(2/\sqrt{d}, 1/2)$-regular subset $B$ of $V(G)$. Since $G$ is $d$-regular and $d > 18\sqrt{d}\log d$, every vertex $u \in V(G)$ satisfies
	\begin{equation}\label{eq:S}
		\sqrt{d} \,\leq\, \deg_B(u) \, \leq\, 3\sqrt{d}.
	\end{equation}
	Set $A \defeq V(G) \setminus B$. We will find a legal removal scheme $(\prec, \mathsf{save})$ on $G$ such that:
	\begin{enumerate}[label=\ep{\itshape\alph*}]
		\item\label{item:RA} In the ordering $\prec$, every vertex in $A$ comes before every vertex in $B$.
		
		\item\label{item:RB} Every vertex in $B$ is saved at least $c \sqrt{d}$ times.
	\end{enumerate}
	Notice that if $(\prec, \mathsf{save})$ satisfies conditions \ref{item:RA} and \ref{item:RB}, then $\mathsf{gap}(\prec, \mathsf{save}) \geq c\sqrt{d}$, which is the property we want. Indeed, take any vertex $u \in V(G)$. If $u \in A$, then, by \ref{item:RA} and \eqref{eq:S},
	\[
	\mathsf{gap}(u; \prec, \mathsf{save}) \,\geq\, |\set{v \in N_G(u) \,:\, v \succ u}| \,\geq\, \deg_B(u) \,\geq\, \sqrt{d}.
	\]
	On the other hand, if $u \in B$, then, by \ref{item:RB},
	\[
	\mathsf{gap}(u; \prec, \mathsf{save}) \,\geq\, |\set{v \in N_G(u) \,:\, \mathsf{save}(v) = u}| \,\geq\, c\sqrt{d}.
	\]
	Assuming $c < 1$, we have $\mathsf{gap}(u; \prec, \mathsf{save}) \geq c\sqrt{d}$ in both cases, as desired.
	
	A legal removal scheme $(\prec, \mathsf{save})$ satisfying \ref{item:RA} and \ref{item:RB} is constructed as follows. For $1 \leq i \leq k$, we recursively define the following numerical parameters:
	\[
	N_1 \,\defeq\, 1 \qquad \text{and} \qquad N_i \,\defeq\, 20k \sum_{j = 1}^{i-1} N_j \text{ for } i \geq 2.
	\]
	Set $p_i \defeq N_i/(6kN_k)$. Note that $p_1 < p_2 < \cdots < p_k = 1/(6k)$. We shall assume $c$ is so small that
	\begin{equation}\label{eq:csmall}
		32 k c \,<\, p_1.
	\end{equation}
	Since $\chi(G) \leq k$, we can partition $A$ into $k$ independent sets $A_1$, \ldots, $A_k$. Let $C_i$ be a $(p_i, 1/2)$-regular subset of $A_i$ and let $D_i$ be a $(p_i, 1/2)$-regular subset of $A_i \setminus C_i$. The ordering $\prec$ is defined by listing the elements of $V(G)$ in the following order:
	\[
	C_1, \ D_1, \ C_2, \ D_2, \ \ldots, \ C_k, \ D_k, \ A \setminus \bigcup_{i=1}^k (C_i \cup D_i), \ B.
	\]
	\ep{The order of the elements in each set in this list is arbitrary.} Since the elements of $B$ appear last in this ordering, condition \ref{item:RA} is fulfilled.
	
	Now we need to define the function $\mathsf{save}$ so that condition \ref{item:RB} holds. We start by recording the following observation:
	
	\begin{claim}\label{claim:CDupper}
		Every vertex $u \in V(G)$ satisfies
		\[
		\deg_{C_i}(u) \,\leq\, \frac{3p_i}{2}d \qquad \text{and} \qquad \deg_{D_i}(u) \,\leq\, \frac{3p_i}{2} d.
		\]
	\end{claim}
	\begin{scproof}
		Immediate from the definitions of $C_i$ and $D_i$ and since the maximum degree of $G$ is $d$. 
	\end{scproof}
	
	By \eqref{eq:S}, each vertex $u \in V(G)$ has at least $d - 3\sqrt{d} \geq d/2$ neighbors in $A$. Therefore, we may partition $B$ into $k$ sets $B_1$, \ldots, $B_k$ so that each vertex in $B_i$ has at least $d/(2k)$ neighbors in $A_i$. This implies that every vertex in $B_i$ has many neighbors in $C_i$ and $D_i$.
	
	\begin{claim}\label{claim:CDlower}
		Every vertex $w \in B_i$ satisfies
		\[
		\deg_{C_i}(w) \,\geq\, \frac{p_i}{4k}d \qquad \text{and} \qquad \deg_{D_i}(w) \,\geq\, \frac{p_i}{8k}d.
		\]
	\end{claim}
	\begin{scproof}
		The first inequality holds since $C_i$ is a $(p_i, 1/2)$-regular subset of $A_i$ and $\deg_{A_i}(w) \geq d/(2k)$. The second inequality follows similarly since, by Claim~\ref{claim:CDupper},
		\[
		\deg_{A_i \setminus C_i}(w) \,\geq\, \frac{d}{2k} - \frac{3p_i}{2}d \,\geq\, \frac{d}{4k}. \qedhere
		\]
	\end{scproof}

	Note that, by \eqref{eq:S}, each vertex in $D_i$ has at most $3 \sqrt{d}$ neighbors in $B_i$. On the other hand, by Claim~\ref{claim:CDlower}, every vertex in $B_i$ has at least $p_id/(8k) \geq p_1 d/(8k)$ neighbors in $D_i$. Since, by \eqref{eq:csmall},
	\[
	\lceil c \sqrt{d} \rceil \cdot 3 \sqrt{d} \,<\, 4cd \, < \, \frac{p_1}{8k} d,
	\]
	we can apply Lemma~\ref{lemma:Hall} to find a partial function $s_i \colon D_i \pto B_i$ such that:
	\begin{itemize}
		\item for all $u \in D_i$, if $s_i(u)$ is defined, then $s_i(u)$ is a neighbor of $u$;
		
		\item the preimage of every vertex $w \in B_i$ under $s_i$ has cardinality $\lceil c \sqrt{d} \rceil$.
	\end{itemize}
	Now we can define $\mathsf{save} \colon V(G) \pto V(G)$ by
	\[
	\mathsf{save}(u) \,\defeq\, \begin{cases}
		s_i(u) &\text{if } u \in D_i \text{ and } s_i(u) \text{ is defined};\\
		\mathsf{blank} &\text{otherwise}.
	\end{cases}
	\]
	By the choice of $s_i$, $(\prec, \mathsf{save})$ is a removal scheme that satisfies \ref{item:RB}. It remains to verify that this removal scheme is legal. To this end, take any $u$, $w \in V(G)$ such that $\mathsf{save}(u) = w$. By construction, this means that $u \in D_i$ and $w \in B_i$ for some $i$. The vertices that precede $u$ in the ordering $\prec$ are the ones in $C_1$, $D_1$, \ldots, $C_{i-1}$, $D_{i-1}$, $C_i$, plus possibly some vertices in $D_i$. Since the set $A_i$ is independent, $u$ has no neighbors in $C_i \cup D_i$, and hence, by Claim~\ref{claim:CDupper},
	\[
	|\set{v \in N_G(u) \,:\, v \prec u}| \,=\, \sum_{j = 1}^{i-1} (\deg_{C_j}(u) + \deg_{D_j}(u)) \,\leq\, \sum_{j = 1}^{i-1}3p_jd \,=\, \frac{3p_i}{20k}d \,<\, \frac{p_i}{4k}d.
	\]
	On the other hand, since no vertex in $C_i$ saves $w$, Claim~\ref{claim:CDlower} yields
	\[
	|\set{v \in N_G(w) \,:\, v \prec u \text{ and } \mathsf{save}(v) \neq w}| \,\geq\, \deg_{C_i}(w) \,\geq\, \frac{p_i}{4k}d.
	\]
	Therefore, inequality \eqref{eq:legal} holds, and the proof of Theorem~\ref{theo:chrom1} is complete.
	
	\subsection{Graphs of girth at least 5}\label{sec:girth5}

	\begin{theocopy}{theo:girth5}\label{theo:girth51}
		There exist $c > 0$ and $d_0 \in \N$ such that if $G$ is a graph of maximum degree $d \geq d_0$ and girth at least $5$, then $\mathsf{wd}(G) \leq d - c \sqrt{d}$.
	\end{theocopy}
	
	Let $G$ be a graph of maximum degree $d$ and girth at least $5$, where $d$ is sufficiently large. Upon replacing $G$ with a supergraph if necessary, we may assume that $G$ is $d$-regular \ep{Lemma~\ref{lemma:embgirth}}. Let $c$ be a sufficiently small positive constant. As in the proof of Theorem~\ref{theo:chrom1}, we will construct a legal removal scheme $(\prec, \mathsf{save})$ on $G$ such that $\mathsf{gap}(\prec, \mathsf{save}) \geq c\sqrt{d}$. By Lemma~\ref{lemma:remove}, this will yield the desired result.

	By Lemma~\ref{lemma:regular}, there is a $(2/\sqrt{d}, 1/2)$-regular subset $B$ of $V(G)$. Then for every vertex $u \in V(G)$,
	\begin{equation}\label{eq:S1}
		\sqrt{d} \,\leq\, \deg_B(u) \,\leq\, 3\sqrt{d}.
	\end{equation}
	Set $A \defeq V(G) \setminus B$. Every vertex in $A$ has at most $3\sqrt{d}$ neighbors in $B$, while every vertex in $B$ has at least $d - 3\sqrt{d} \geq d/2$ neighbors in $A$. Since $\lceil \sqrt{d}/8\rceil \cdot 3\sqrt{d} < d/2$, Lemma~\ref{lemma:Hall} gives a partial function $s \colon A \pto B$ such that:
	\begin{itemize}
		\item for all $u \in A$, if $s(u)$ is defined, then $s(u)$ is a neighbor of $u$;
		\item the preimage of each $w \in B$ under $s$ has cardinality $\lceil \sqrt{d}/8\rceil$.
	\end{itemize}
	For each $w \in B$, let $S_w$ denote the preimage of $w$ under $s$; for $u \in A$, set $S_u \defeq \0$.
	
	Now we assemble a removal scheme $(\prec, \mathsf{save})$ using a randomized procedure. Pick a linear ordering $\lhd$ of $A$ uniformly at random. The ordering $\prec$ will start with the vertices of $A$ listed according to $\lhd$, followed by the vertices of $B$ in some order \ep{to be specified shortly}. Intuitively, we imagine that every vertex $u \in A$ attempts to save the vertex $s(u)$. This attempt only succeeds if condition \eqref{eq:legal} is satisfied. Formally, we say that $u \in A$ with $s(u) = w \in B$ is \emphd{successful} if
	\[
	|\set{v \in N_A(u) \,:\, v \lhd u}| \,<\, |\set{v \in N_A(w) \setminus S_w \,:\, v \lhd u}|.
	\]
	If $u \in A$ is successful, then we write $\mathsf{save}(u) \defeq s(u)$; for all other vertices $u$ we set $\mathsf{save}(u) \defeq \mathsf{blank}$. 
	
	Say that a vertex $w \in B$ is \emphd{happy} if its preimage under the function $\mathsf{save}$ has cardinality at least $c \sqrt{d}$. Let $H \subseteq B$ be the set of all happy vertices. The ordering $\prec$ consists of $A$ listed according to $\lhd$, followed by $B \setminus H$ in an arbitrary order, and then by $H$ in an arbitrary order. By construction, $(\prec, \mathsf{save})$ is a legal removal scheme, and we claim that $\mathsf{gap}(\prec,\mathsf{save}) \geq c \sqrt{d}$ with positive probability. The key fact we need to establish is the following:
	
	\begin{claim}\label{claim:happy}
		With positive probability, every vertex of $G$ has at least $c\sqrt{d}$ neighbors in $H$.
	\end{claim}
	
	Let us see why Claim~\ref{claim:happy} implies the desired result. Suppose that every vertex of $G$ has at least $c\sqrt{d}$ neighbors in $H$. Take any $u \in V(G)$. If $u \not \in H$, then
	\[
	\mathsf{gap}(u; \prec, \mathsf{save}) \,\geq\, |\set{v \in N_G(u) \,:\, v \succ u}| \,\geq\, \deg_H(u) \,\geq\, c\sqrt{d}.
	\]
	On the other hand, if $u \in H$, then, by the definition of $H$,
	\[
	\mathsf{gap}(u; \prec, \mathsf{save}) \,\geq\, |\set{v \in N_G(u) \,:\, \mathsf{save}(v) = u}| \,\geq\, c\sqrt{d}.
	\]
	In either case, $\mathsf{gap}(u; \prec, \mathsf{save}) \geq c\sqrt{d}$, as desired.
	
	In the remainder of this section we prove Claim~\ref{claim:happy}. It will be convenient to assume that the random ordering $\lhd$ is sampled according to the following procedure: each vertex $u \in A$ picks a real number $\theta(u) \in [0,1]$ uniformly at random, and then we set $u_1 \lhd u_2$ if and only if $\theta(u_1) < \theta(u_2)$ \ep{note that $\theta(u_1) \neq \theta(u_2)$ with probability $1$}. For each $u \in V(G)$, let $X_u$ be the random event that $\deg_H(u) < c \sqrt{d}$. It is clear that $X_u$ only depends on the values of the function $\theta$ on the vertices at distance at most $3$ from $u$. Therefore, $X_u$ is mutually independent from the events $X_v$ corresponding to the vertices $v$ at distance more than $6$ from $u$. Hence, by the \hyperref[theo:LLL]{Lov\'asz Local Lemma}, to prove that with positive probability none of the events $X_u$ happen it suffices to show that
	\begin{equation}\label{eq:6}
		\P[X_u] \,=\, o(d^{-6}).
	\end{equation}
	
	The proof of \eqref{eq:6} is somewhat technical, so before getting into its details, let us briefly explain the intuition behind our approach. Assuming $c$ is small enough, it is possible to show that for each $w \in B$, $\P[\text{$w$ is happy}] = \Omega(1)$. Since every vertex $u \in V(G)$ has at least $\sqrt{d}$ neighbors in $B$, we have $\E[\deg_H(u)] = \Omega(\sqrt{d})$. Ideally, we would now argue that the random variable $\deg_H(u)$ is close to its expected value with very high probability. One way to achieve this would be to show that the random events ``$w$ is happy'' for $w \in N_B(u)$ are close to being mutually independent and then apply the Chernoff bound or some other similar result. This strategy indeed works in the case when $G$ has girth at least $7$. This is because for each $w \in N_B(u)$, the event ``$w$ is happy'' is determined by the values of $\theta$ in the radius-$2$ ball around $w$, and the girth-$7$ assumption implies that the radius-$2$ balls around the vertices in $N_B(u)$ do not overlap too much.
	
	It turns out that, with a more clever argument, we can reduce the girth requirement from $7$ to $5$. The idea is to define a certain property of vertices $w \in B$, which we call being \emph{powerful} \ep{or, more accurately, \emph{$\epsilon$-powerful} for some $\epsilon > 0$}, so that the following statements hold:
	\begin{enumerate}[label=\ep{\itshape\alph*}]
		\item\label{item:local} the event ``$w$ is powerful'' is determined by the values of $\theta$ on the neighbors of $w$;
		
		\item\label{item:constant} the probability that $w$ is powerful is at least $\Omega(1)$ \ep{Claim~\ref{claim:strong}};
		
		\item\label{item:implieshappy} if $w$ is powerful, then $w$ is happy with very high probability \ep{Claim~\ref{claim:stronguse}}.
	\end{enumerate}
	Thanks to \ref{item:constant}, the expected number of powerful neighbors for each vertex $u \in V(G)$ is $\Omega(\sqrt{d})$. Using \ref{item:local} and the girth-$5$ assumption, we can show that in fact $u$ has $\Omega(\sqrt{d})$ powerful neighbors with very high probability. Finally, according to \ref{item:implieshappy}, once $u$ has $\Omega(\sqrt{d})$ powerful neighbors, it is extremely likely that it has $\Omega(\sqrt{d})$ happy neighbors as well.
	
	Let us now begin the formal proof. We start by associating to each vertex of $G$ a \ep{random} vector with entries in $[0,1]$ by setting, for every $u \in V(G)$,
	\[
	x_u \,\defeq\, (\theta(v) \,:\, v \in N_A(u) \setminus S_u).
	\]
	\ep{Recall that $S_u = \0$ for $u \in A$.} Now we introduce the following definitions:
	
	\begin{defn}[\textls{Powerful vectors and vertices}]
		Given a vector $x = (x_1, \ldots, x_k) \in [0,1]^k$ and a real number $\alpha \in [0,1]$, let the \emphd{$\alpha$-power} of $x$ be the quantity
		\[
		\pi(x, \alpha) \,\defeq\, |\set{i \,:\, x_i < \alpha}|.
		\]
		For $\epsilon > 0$, we say that a vector $x \in [0,1]^k$ is \emphd{$\epsilon$-powerful} if the following statement holds: If we pick a real number $\alpha \in [0,1]$ and a vector $y \in [0,1]^d$ uniformly at random, then
		\begin{equation}\label{eq:power}
			\P[\pi(y, \alpha) < \pi(x, \alpha)] \,\geq\, \epsilon.
		\end{equation}
		A vertex $w \in B$ is \emphd{$\epsilon$-powerful} if the vector $x_w$ is $\epsilon$-powerful.
	\end{defn}
	
	We remark that if $x \in [0,1]^k$ is $\epsilon$-powerful, then \eqref{eq:power} also holds for $y$ drawn uniformly at random from $[0,1]^\ell$ for any $\ell \leq d$. Similarly, if an $\epsilon$-powerful vector $x$ is obtained from another vector $x'$ by removing some of the coordinates, then $x'$ is $\epsilon$-powerful as well, since $\pi(x',\alpha) \geq \pi(x, \alpha)$ for all $\alpha$.
	
	Using this notation, we can say that a vertex $u \in A$ with $s(u) = w$ is successful if and only if
	\[
	\pi(x_u, \theta(u)) \,<\, \pi(x_w, \theta(u)).
	\]

	\begin{claim}\label{claim:strong}
		There exists a constant $\epsilon > 0$ such that if $k \geq d - 5\sqrt{d}$, then the probability that a uniformly random vector $x \in [0,1]^k$ is $\epsilon$-powerful is at least $\epsilon$.
	\end{claim}
	\begin{scproof}
		For $\epsilon > 0$, let $p_\epsilon$ denote the probability that a uniformly random vector $x \in [0,1]^k$ is $\epsilon$-powerful. If we sample $x \in [0,1]^k$, $\alpha \in [0,1]$, and $y \in [0,1]^d$ uniformly at random, then
		\begin{align}
			\P[\pi(y, \alpha) < \pi(x, \alpha)] \,=\, & \P[\text{$x$ is $\epsilon$-powerful}] \P[\pi(y, \alpha) < \pi(x, \alpha) \,\vert\, \text{$x$ is $\epsilon$-powerful}] \nonumber\\
			&\,+\, \P[\text{$x$ is not $\epsilon$-powerful}] \P[\pi(y, \alpha) < \pi(x, \alpha) \,\vert\, \text{$x$ is not $\epsilon$-powerful}] \nonumber\\
			\leq\, & p_\epsilon + (1-p_\epsilon) \epsilon \,\leq\, p_\epsilon + \epsilon. \label{eq:power1}
		\end{align}
		We now prove a lower bound on the left-hand side of \eqref{eq:power1}. We sample $\alpha \in [0,1]$ first. Note that with probability $1/3$, we get $1/3 \leq \alpha \leq 2/3$. Now $\pi(x, \alpha)$ and $\pi(y,\alpha)$ are independent random variables sampled from the binomial distributions $\mathrm{Bin}(k, \alpha)$ and $\mathrm{Bin}(d, \alpha)$ respectively. It follows from the \hyperref[theo:BE]{Berry--Esseen theorem} \ep{specifically from equation \eqref{eq:BE}} that there exists a constant $\gamma > 0$ such that, assuming $d$ is large enough and $1/3 \leq \alpha \leq 2/3$, we have
		\[
		\P[\pi(x, \alpha) > \alpha k] \,\geq\,\gamma \quad \text{and} \quad \P[\pi(y, \alpha) < \alpha (d - 5 \sqrt{d})] \,\geq\, \gamma.
		\]
		Since $\alpha (d - 5\sqrt{d}) \leq \alpha k$, we conclude that
		\[
		\P[\pi(y, \alpha) < \pi(x, \alpha)] \,\geq\, \frac{\gamma^2}{3}. 
		\]
		By \eqref{eq:power1}, setting $\epsilon \defeq \gamma^2/6$ finishes the proof.
	\end{scproof}
	
	In the remainder of the proof we fix a constant $\epsilon$ satisfying the conclusion of Claim~\ref{claim:strong}. We shall assume that the ratio $c/\epsilon$ is sufficiently small, say $c < \epsilon/10$. To simplify the presentation, we will use the asymptotic notation $O(\cdot)$ to hide positive constant factors \ep{which may be computed as functions of $\epsilon$ and $c$}.
	
	\begin{claim}\label{claim:stronguse}
		For every vertex $w \in B$, we have
		\[
		\P\left[\text{$w$ is happy}\,\middle\vert\, \text{$w$ is $\epsilon$-powerful} \right] \,\geq\, 1 - \exp\left(-O(\sqrt{d})\right).
		\]
	\end{claim}
	\begin{scproof}
		Let us fix the values $\theta(v)$ for $v \in N_A(w) \setminus S_w$ so that the vector $x_w$ is $\epsilon$-powerful. Now consider any $u \in S_w$. The value $\theta(u)$ is chosen uniformly at random from $[0,1]$. Moreover, since $G$ is triangle-free, $u$ and $w$ have no common neighbors, which means that the values $\theta(v)$ for $v \in N_A(u)$ have not yet been determined. In other words, $x_u$ is a uniformly random vector from $[0,1]^{\deg_A(u)}$. Since $x_w$ is $\epsilon$-powerful and $\deg_A(u) \leq d$,
		\[
		\P[\text{$u$ is successful}] \,=\, \P[\pi(x_u, \theta(u)) < \pi(x_w, \theta(u))] \,\geq\, \epsilon.
		\]
		Since $G$ has girth at least $5$, the vertices in $S_w$ have no common neighbors except $w$, and thus the random events ``$\pi(x_u, \theta(u)) < \pi(x_w, \theta(u))$'' for $u \in S_w$ are mutually independent. Therefore, the random variable $\xi$ equal to the cardinality of the preimage of $w$ under the function $\mathsf{save}$ is bounded below by a binomial random variable with distribution $\mathrm{Bin}(|S_w|, \epsilon)$. Hence, we may apply the \hyperref[theo:Chernoff]{Chernoff bound} \ep{Theorem~\ref{theo:Chernoff}} and the inequality $|S_w| \geq \sqrt{d}/8$ to conclude that
		\[
		\P[\text{$w$ is not happy}] \,=\, \P[\xi < c\sqrt{d}] \,<\, 2\exp\left(-\left(\frac{\epsilon}{8}-c\right)^2 \frac{8\sqrt{d}}{3\epsilon}\right) \,\leq\, \exp\left(-O(\sqrt{d})\right). \qedhere
		\]
	\end{scproof}
	
	For a vertex $u \in V(G)$, define \[P_\epsilon(u) \,\defeq\, \set{w \in N_B(u) \,:\, \text{$w$ is $\epsilon$-powerful}}.\]
	
	\begin{claim}\label{claim:P}
		For every vertex $u \in V(G)$, we have
		\[
		\P[|P_\epsilon(u)| \geq c \sqrt{d} ] \,\geq\, 1 - \exp\left(-O(\sqrt{d})\right).
		\]
	\end{claim}
	\begin{scproof}
		A slight technical issue here arises from the fact that the vectors $x_w$ for $w \in N_B(u)$ may not be probabilistically independent from each other, since each of them may include $\theta(u)$ as one of the coordinates. To remedy this, we define for every $w \in N_B(u)$ a vector $x_w'$ as follows:
		\[
		x_w' \,\defeq\, (\theta(v) \,:\, v \in N_A(w) \setminus (S_w \cup \set{u})).
		\]
		That is, $x_w'$ is obtained from $x_w$ by deleting the coordinate corresponding to $u$. Let \[P_\epsilon'(u) \,\defeq\, \set{w \in N_B(u) \,:\, \text{$x_w'$ is $\epsilon$-powerful}}.\] Then $P_\epsilon'(u) \subseteq P_\epsilon(u)$, so it suffices to argue that
		\[
		\P[|P_\epsilon'(u)| \geq c\sqrt{d}] \,\geq\, 1 - \exp\left(-O(\sqrt{d})\right).
		\]
		For $w \in N_B(u)$, let $k(w) \defeq |N_A(w) \setminus (S_w \cup \set{u})|$. Then, by \eqref{eq:S1} and since $|S_w| = \lceil \sqrt{d}/8 \rceil$, we have
		\[
		k(w) \,\geq\, d \,-\, 3\sqrt{d} \,-\, \lceil \sqrt{d}/8\rceil \,-\, 1 \,\geq\, d - 5\sqrt{d}.
		\]
		By the choice of $\epsilon$ and since $x_w'$ is drawn uniformly at random from $[0,1]^{k(w)}$, we conclude that
		\[
		\P[\text{$x_w'$ is $\epsilon$-powerful}] \,\geq\, \epsilon.
		\]
		As $G$ has girth at least $5$, the vertices in $N_B(u)$ have no common neighbors except $u$, so we can apply the \hyperref[theo:Chernoff]{Chernoff bound} and the inequality $\deg_B(u) \geq \sqrt{d}$ to obtain the desired bound
		\[
		\P[|P_\epsilon'(u)| < c\sqrt{d}] \,<\, 2\exp\left(- \frac{(\epsilon - c)^2\sqrt{d}}{3\epsilon}\right) \,\leq \, \exp\left(-O(\sqrt{d})\right). \qedhere
		\]
	\end{scproof}
	
	Finally, we can bound the probability of each event $X_u$:
	
	\begin{claim}\label{claim:final}
		Let $u \in V(G)$. Recall that $X_u$ is the event that $\deg_H(u) < c \sqrt{d}$. Then \[\P[X_u] \,\leq\, \exp\left(-O(\sqrt{d})\right).\]
	\end{claim}
	\begin{scproof}
		By Claim~\ref{claim:stronguse}, for each $w \in N_B(u)$, we have
		\[
		\P[w \in P_\epsilon(u) \setminus H] \,=\, \P\left[\text{$w$ is $\epsilon$-powerful but not happy} \right] \,\leq\, \exp\left(-O(\sqrt{d})\right).
		\]
		Therefore, by the union bound,
		\[
		\P[P_\epsilon(u) \setminus H \neq \0] \,\leq\, \deg_B(u) \cdot  \exp\left(-O(\sqrt{d})\right) \,\leq\, \exp\left(-O(\sqrt{d})\right).
		\]
		Hence, by Claim~\ref{claim:P},
		\[
		\P[X_u] \, \leq \, \P[|P_\epsilon(u)| < c \sqrt{d} ] \,+\, \P[P_\epsilon(u) \setminus H \neq \0] \,\leq\, \exp\left(-O(\sqrt{d})\right). \qedhere
		\]
	\end{scproof}
	
	The upper bound on $\P[X_u]$ given by Claim~\ref{claim:final} implies the asymptotic bound \eqref{eq:6}. As discussed earlier, this yields Claim~\ref{claim:happy} and completes the proof of Theorem~\ref{theo:girth51}.
	
	\printbibliography

\end{document}